# Probability & incompressible Navier-Stokes equations: An overview of some recent developments[*]

## Edward C. Waymire[†]

*Department of Mathematics, Oregon State University, Corvallis, OR 97370, USA e-mail:*
`waymire@math.orst.edu`

**Abstract:** This is largely an attempt to provide probabilists some orientation to an important class of non-linear partial differential equations in applied mathematics, the incompressible Navier-Stokes equations. Particular focus is given to the probabilistic framework introduced by LeJan and Sznitman [Probab. Theory Related Fields 109 (1997) 343–366] and extended by Bhattacharya et al. [Trans. Amer. Math. Soc. 355 (2003) 5003–5040; IMA Vol. Math. Appl., vol. 140, 2004, in press]. In particular this is an effort to provide some foundational facts about these equations and an overview of some recent results with an indication of some new directions for probabilistic consideration.

## 1. Navier-Stokes Equations: Background

The present paper is an attempt to provide a brief orientation to Navier-Stokes equations from a probabilistic perspective developed in the course of working with a focussed research group in this area during the past few years at Oregon State University (OSU). An effort is made to furnish a bit more background for the uninitiated on some of the basics of these equations along with a summary description of the probabilistic methods and results. The approach is in somewhat "broad strokes". The reader should be able to follow and/or supply most basic calculations, but the more technical proofs of some of the main results are left to original references. More detailed accounts of some of the mathematical intricacies, conjectures and open problems identified by the OSU group can be found in Bhattacharya et al. [6, 7], Chen et al. [12], and in a recently completed PhD thesis by Orum [45].

---

[*]This research was partially supported by Focussed Research Group collaborative award DMS-0073958 and by DMS-0327705 to Oregon State University. This article was completed during a sabbatical year (2004-2005) visiting appointment in the Department of Mathematics at Cornell University.

[†]This article is based on collaboration with my Oregon State University colleagues Larry Chen, Ronald Guenther, Mina Ossiander, Enrique Thomann, OSU PhD student Chris Orum, now at Los Alamos National Laboratory, and Rabi Bhattacharya at the University of Arizona and Indiana University. Any new insights provided in this article are the direct result of ongoing close collaboration with this research group.

A purely mathematical description of the Navier-Stokes equations naturally depends on the framework (function space) in which initial data and external forcings are provided, and in which solutions are sought, e.g. see see Ladyzhenskaya [32, 33], Temam [52]. However, rather than begin with the technical nomenclature, let us first take a look at the underlying physical picture. The 3-dimensional incompressible Navier-Stokes (NS) equations governing fluid velocities $\mathbf{u}(\mathbf{x}, t) = (u_1(\mathbf{x}, t), u_2(\mathbf{x}, t), u_3(\mathbf{x}, t))$, $\mathbf{x} \in \mathbf{R}^3$, $t \geq 0$, with initial data $\mathbf{u}(\mathbf{x}, 0^+) = \mathbf{u}_0(\mathbf{x})$, are a manifestation of Newton's law of motion and mass conservation for a homogeneous (constant density) fluid:

$$\frac{\partial \mathbf{u}}{\partial t} + (\mathbf{u} \cdot \nabla)\mathbf{u} = \nu \Delta \mathbf{u} - \nabla p + \mathbf{g}, \quad \nabla \cdot \mathbf{u} = 0, \quad \mathbf{u}(\mathbf{x}, 0^+) = \mathbf{u}_0(\mathbf{x}), \quad (1)$$

where $\nabla = (\partial/\partial x_1, \partial/\partial x_2, \partial/\partial x_3)$, $\mathbf{u} \cdot \nabla = \sum_{j=1}^{3} u_j \frac{\partial}{\partial x_j}$, and $\Delta = \nabla \cdot \nabla = \sum_{j=1}^{3} \partial^2/\partial x_j^2$; the latter two are applied component-wise to vector-valued functions. The term $p(\mathbf{x}, t)$ is the (scalar) *pressure*, $\mathbf{g}(\mathbf{x}, t)$ represents external *forcing*, and $\nu > 0$ is the kinematic *viscosity*. Probabilists might want to write the term $\nu \Delta$ as $\frac{1}{2} 2\nu \Delta$ in anticipation of the *heat kernel* in the form $k(\mathbf{x}, t) = (4\pi\nu t)^{-\frac{3}{2}} \exp\{-\frac{|\mathbf{x}|^2}{4\nu t}\}$ more familiar in applied mathematics. The *incompressibility* condition refers to the requirement that

$$\nabla \cdot \mathbf{u} = 0. \quad (2)$$

The (in-viscid) equations corresponding to $\nu = 0$ are referred to as the *Euler equations,* and are generally considered apart from Navier-Stokes.

The equations (1) comprise a (non-linear) system of four partial differential equations governing four unknown scalar quantities; three coordinates of velocity $u_1, u_2, u_3$ and a scalar pressure $p$. The equations model the motion of a viscous incompressible homogeneous fluid in an inertial frame in free space. They are derived from conservation of mass and momenta of "fluid elements" (in the sense of continuum mechanics) under the additional assumption that there is a local linear relation between the stress and strain rates experienced by fluid parcels. In (1), the left-hand side $\frac{\partial \mathbf{u}}{\partial t} + (\mathbf{u} \cdot \nabla)\mathbf{u}$ represents the acceleration in an Eulerian frame of reference; i.e. one views the motion as that of a "parade of fluid elements" and $\mathbf{u}(\mathbf{x}, t)$ is the velocity of the fluid parcel at that location and at that time. Thus if a parcel has moved to location $\mathbf{x}(t) = (x(t), y(t), z(t))$ at time $t$, then the velocity of the parcel is given by $\mathbf{u}(\mathbf{x}, t) \equiv (u_1(x, y, z, t), u_2(x, y, z, t), u_3(x, y, z, t)) = (\frac{dx}{dt}, \frac{dy}{dt}, \frac{dz}{dt})$ and, using the chain rule, the acceleration is therefore given by the so-called *material* or *convective derivative*

$$\frac{d\mathbf{u}}{dt} = \frac{\partial \mathbf{u}}{\partial t} + \frac{\partial \mathbf{u}}{\partial x}\frac{dx}{dt} + \frac{\partial \mathbf{u}}{\partial y}\frac{dy}{dt} + \frac{\partial \mathbf{u}}{\partial z}\frac{dz}{dt} = \frac{\partial \mathbf{u}}{\partial t} + (\mathbf{u} \cdot \nabla)\mathbf{u}. \quad (3)$$

The right side of (1) resolves the forces acting on the fluid parcel. The first term $\nu\Delta\mathbf{u}$ represents the viscous friction forces and is often referred to as the *dissipative term*. This term arises in the fluid mechanical derivation as the result of the aforementioned ad-hoc model of linearized shear stresses due to interactions between molecules comprising the fluid. The term $\mathbf{g} \equiv \mathbf{g}(x, y, z, t)$ represents body forces, such as gravity and/or other externally applied force fields. The pressure gradient term $-\nabla p$ defines the local pressure force on a fluid parcel (depending on $p$ up to an arbitrary constant). Physically the pressure gradient is a force which acts to maintain a homogeneous fluid density under incompressibility. Thus one may reasonably expect that pressure can be computed from the velocity and boundary conditions (at infinity in this case). We will see that the following orthogonality opens a way to mathematically first remove the pressure for focus on the incompressible velocity. Integrating the following "product formula" over a ball $B_0(R)$ of radius $R$

$$\nabla \cdot (p\mathbf{u}) = \mathbf{u} \cdot \nabla p + (\nabla \cdot \mathbf{u})p \tag{4}$$

and then applying Gauss' divergence theorem yields the following version of an "integration by parts" type formula

$$\int_{B_0(R)} \mathbf{u} \cdot \nabla p = \int_{\partial B_0(R)} p\mathbf{u} \cdot \mathbf{n} - \int_{B_0(R)} (\nabla \cdot \mathbf{u})p. \tag{5}$$

Thus, if $\mathbf{u}$ is incompressible then the second integral on the right side of (5) vanishes and the first integral is $O(\|p\mathbf{u}\mathbf{1}(B_0(R))\|_\infty R^2)$; here $\mathbf{1}(B_0(R))$ denotes the indicator function. Assuming $p\mathbf{u}$ decays sufficiently fast at infinity, one therefore obtains the following basic $L^2-$orthogonality

$$\int_{\mathbf{R}^3} \mathbf{u} \cdot \nabla p = 0. \tag{6}$$

**Remark**. Decay assumptions are implicit to most of the theory to be described throughout this article, from orthogonality conditions (6) to Fourier transform formulations. Such boundary conditions at infinity on Navier-Stokes are the subject of extensive investigation for both free-space and exterior domain problems since Finn's conception of the notion of "physically reasonable" solution; e.g. see Finn [18], Heywood [28], Galdi [23]. A perhaps more familiar illustration of the need to impose decay conditions occurs for example in the simple case of the heat equation, where one observes that since the Gaussian heat kernel is a tempered function, its convolution with a tempered distribution is well-defined as a function. In particular, it is sufficient to assume that the initial data is a tempered distribution for the well-known solution as a convolution with the heat kernel to apply.

In summary, two distinguishing features of the equations (1) are: (a) the nonlinear term, referred to as the *convective term*, $(\mathbf{u} \cdot \nabla)\mathbf{u}$ on the left-hand side, whose role in representing the acceleration is intrinsic and cannot be dropped by changing the model of the forces on the right-hand side, and (b) the incompressibility condition $\nabla \cdot \mathbf{u} = 0$. This latter condition is an additional equation to momentum balance obtained by consideration of mass conservation for a homogeneous fluid of constant fluid density.

Systematic derivations of the Navier-Stokes equations can be found in any number of standard texts on fluids, e.g. Batchelor [3], Landau and Lifschitz [34], Chorin and Marsden [15]. Thumbnail historical sketches of Claude Louis Marie Henri Navier (1785-1836) and George Gabriel Stokes (1819-1903) are available on the web. A brief textual history as it pertains to fluid flow as well as a fundamentally significant role in engineering design can be found in Anderson [1].

Let us briefly turn to the mathematical framework in which one might view these equations. Cannone [9] contains an excellent survey of a variety of natural function spaces in which to analyze the equations. However, rather than attempt to develop the precise function space requirements at the outset for the results to be discussed, we will for the most part continue to confine ourselves to $L^2(\mathbf{R}^3)^3$ and/or spaces of *solenoidal*, i.e. $\nabla \cdot \mathbf{u} = 0$, tempered distributions $\mathbf{u}$, where the reader will be able to check basic statements and calculations using standard Fourier methods as covered, for example, in Folland [19]. The actual function spaces will evolve naturally in the course of the development of the probabilistic approach.

An in-depth analysis of these equations given by Jean Leray [36][1] has become a lasting benchmark for contemporary research. Under the assumptions that $\mathbf{u}_0, \mathbf{g} \in L^2(\mathbf{R}^3)^3$, $\nabla \cdot \mathbf{u}_0 = 0$, and an *energy inequality* of the form

$$1/2||\mathbf{u}(\cdot,t)||_2^2 + \int_0^t ||\nabla \cdot \mathbf{u}(\cdot,s)||_2^2 ds \leq 1/2||\mathbf{u}(\cdot,0)||_2^2 + \int_0^t < \mathbf{g}(\cdot,s), \mathbf{u}(\cdot,s) > ds \tag{7}$$

applies (in the distributional sense), Leray [36] proves the existence of a global (weak) solution to (1) satisfying the inequality (7) . The corresponding uniqueness in $L^2(\mathbf{R}^3)^3$ of Leray's solution remains open to this day for the case of $n = 3$ space dimensions. From physical considerations, namely conservation of energy, the energy inequality (7) is expected to be an equality. However it has not been proven that such a solution exists in three-dimensions. A formal derivation of conservation of energy is obtained by taking the scalar product of the terms in (1) with $\mathbf{u}$ and then integrating over $\mathbf{R}^3$. Assuming sufficiently rapid decay of $\mathbf{u}, p$ at infinity, integration by parts and incompressibility yield that the convective term and pressure both vanish. For the scalar product with the dissipative term one has similarly by performing integration by parts

$$\nu \int_{\mathbf{R}^3} \Delta \mathbf{u} \cdot \mathbf{u} = \nu \sum_{i,j=1}^3 \int_{\mathbf{R}^3} \frac{\partial^2 u_i}{\partial x_j^2} u_j = -\nu \sum_{i,j=1}^3 \int_{\mathbf{R}^3} (\frac{\partial u_i}{\partial x_j})^2 = -\nu ||\nabla \cdot \mathbf{u}(\cdot,t)||_2^2. \tag{8}$$

Noting that $\frac{1}{2}\frac{d}{dt}\int_{\mathbf{R}^3} |\mathbf{u}(\cdot,t)|^2 \equiv \int \frac{\partial \mathbf{u}}{\partial t} \cdot \mathbf{u}$ one has the equation for conservation of energy

$$\frac{1}{2}\frac{d}{dt}||\mathbf{u}(\cdot,t)||_2^2 + \nu||\nabla \cdot \mathbf{u}(\cdot,t)||_2^2 = < \mathbf{g}, \mathbf{u} > . \tag{9}$$

---

[1]Robert Terrell has made his translation of Leray's seminal paper to English available on the web at http://www.math.cornell.edu/~bterrell/leray.shtml

At another extreme for solution spaces, in defining the Millenium Prize Problem for the Clay Institute, Fefferman [17] stipulates that a physically reasonable solution must be infinitely smooth and decrease more rapidly than any power of the spatial variable $\mathbf{x}$. However this choice of solution space is not without controversy in the international mathematical community. In a critique of this requirement, Ladyszhenskya [33] remarks that "...one must leave both the choice of the phase space and the class of generalized solutions to the researcher without prescribing to him infinite smoothness or some other smoothness of solutions. The only requirement needed is indeed that the uniqueness theorem must hold in the chosen class of generalized solutions." Here she is addressing her basic point that the first fundamental question is whether the equations (and boundary conditions) do uniquely determine the dynamics of an incompressible fluid. In particular, questions regarding the qualitative behavior and properties of the solution naturally follow. In the context of this article boundary conditions refer to decay properties at infinity, however the methods extend to periodic boundary as well.[2]

Apart from definitions of classical and various forms of weak (in distributional sense) solutions, some of the mathematical intricacies are further illustrated with the work of Fujita and Kato [21], Kato [30] in which a notion of "mild solution" is advanced as follows: The so-called *Leray projection operator* $\mathcal{P}$ is most easily defined by the orthogonal projection of the Hilbert space $L^2(\mathbf{R}^3)^3$ onto the closed subspace $\{\mathbf{u} \in L^2(\mathbf{R}^3)^3 : \nabla \cdot \mathbf{u} = 0\}$. In particular, denoting (componentwise) Fourier transforms on $\mathbf{u} = (u_1, u_2, u_3) \in L^2(\mathbf{R}^3)^3$ by $\hat{\mathbf{u}} = (\hat{u}_1, \hat{u}_2, \hat{u}_3)$, where for $f \in L^2(\mathbf{R}^3)$

$$\hat{f}(\xi) = \frac{1}{(2\pi)^{\frac{3}{2}}} \int_{\mathbf{R}^3} e^{-i\mathbf{x}\cdot\xi} f(\mathbf{x}) d\mathbf{x}, \tag{10}$$

the subspace of divergence-free (incompressible) vector fields corresponds to the following basic orthogonality in Fourier frequency space

$$\xi \cdot \widehat{\mathcal{P}\mathbf{u}}(\xi) = 0. \tag{11}$$

Now, since orthogonality of the projection requires that $\hat{\mathbf{u}}(\xi) - \widehat{\mathcal{P}\mathbf{u}}(\xi)$ be orthogonal to $\widehat{\mathcal{P}\mathbf{u}}(\xi)$, this vector difference is in the direction of $\xi$, i.e. $\hat{\mathbf{u}}(\xi) - \widehat{\mathcal{P}\mathbf{u}}(\xi) = (\frac{\xi}{|\xi|} \cdot \hat{\mathbf{u}}(\xi))\frac{\xi}{|\xi|}$. Thus,

$$\widehat{\mathcal{P}\mathbf{u}}(\xi) = \hat{\mathbf{u}}(\xi) - (\frac{\xi}{|\xi|} \cdot \hat{\mathbf{u}}(\xi))\frac{\xi}{|\xi|} = (\mathbf{I} - \frac{\xi \otimes \xi}{|\xi|^2})\hat{\mathbf{u}}(\xi), \quad \xi \in \mathbf{R}^3, \tag{12}$$

where $\mathbf{I} = ((\delta_{ij}))_{1 \le i,j \le 3}$ is the identity matrix and the tensor product $\otimes$ of two vectors $\mathbf{a}, \mathbf{b} \in \mathbf{R}^3$ is defined by the matrix $\mathbf{a} \otimes \mathbf{b} = ((a_i b_j))_{1 \le i,j \le 3}$. The projection extends to other function spaces for which Fourier transforms can be defined. In view of the incompressibility condition and (6) one may apply the projection operator to (1) to obtain for incompressible initial data $\mathbf{u}_0$ belonging to a suitable function space, e.g. $L^2(\mathbf{R}^3)^3$, the "projected equation"

$$\frac{\partial \mathbf{u}}{\partial t} + \mathcal{P}(\mathbf{u} \cdot \nabla \mathbf{u}) = \nu \Delta \mathbf{u} + \mathcal{P}(\mathbf{g}), \quad \mathbf{u}(\mathbf{x}, 0^+) = \mathbf{u}_0(\mathbf{x}). \tag{13}$$

---

[2]For the most part analysts do not make a significant distinction in the choice between free-space and periodic boundary conditions. While both free-space and periodic boundary are amenable to stochastic cascade representations, there are subtle differences due to the role of small non-zero frequencies which do not arise in the wave-number lattice associated with Fourier series.

If one views (13) as a (vector) heat equation forced by the (unknown) term $\mathbf{w} = \mathcal{P}(\mathbf{g} - \mathbf{u} \cdot \nabla \mathbf{u})$ then a standard Duhamel argument, e.g. see Bhattacharya and Waymire [4], may be applied to deduce the following integral equation

$$\mathbf{u}(\mathbf{x}, t) = T_t \mathbf{u}_0(\mathbf{x}) + \int_0^t T_{t-s} \mathbf{w}(\mathbf{x}, s) ds, \tag{14}$$

where $\{T_t \equiv e^{t\nu\Delta} : t \geq 0\}$ is the (self-adjoint) semi-group associated with 3-dimensional Brownian motion $\mathbf{Z} = \{\mathbf{Z}(t) : t \geq 0\}$ with zero drift and diffusion coefficient $2\nu$, i.e. $T_t g(\mathbf{x}) = E_{\mathbf{x}} g(\mathbf{Z}(t)) = \int_{\mathbf{R}^3} g(\mathbf{y}) \frac{1}{(4\pi\nu t)^{\frac{3}{2}}} \exp\{-\frac{1}{4\pi\nu t}|\mathbf{y} - \mathbf{x}|^2\} d\mathbf{y}$, and applied component-wise. Solutions to (14) satisfying $\nabla \cdot \mathbf{u} = 0$ are referred to as *mild solutions*. This is the starting point of Kato's approach in which he proved existence of mild solutions for "small" initial data in the special function space $L^3(\mathbf{R}^3)^3$. Specifically, ignoring external forcing for simplicity, Kato [30] showed that there is a number $A > 0$ such that for $\|u_0\|_{L^3} \leq A$, one has a solution such that $t \to \mathbf{u}(\cdot, t) \in L^3(\mathbf{R}^3)^3$ is continuous and satisfies (14) on every compact time interval $[0, T]$, $T > 0$. However it wasn't until very recently that uniqueness was proven by Furioli, Lemarié-Rieusett and Terraneo [22].

The Picard iteration scheme provides a useful technique naturally associated with the (projected) Navier-Stokes equation. Observe that using incompressibility the non-linear term may be expressed, analogously to $uu_x = \frac{1}{2}(u^2)_x$ in one-dimension, in the equivalent form

$$(\mathbf{u} \cdot \nabla)\mathbf{u} = \nabla \cdot (\mathbf{u} \otimes \mathbf{u}). \tag{15}$$

Then the iteration is given by

$$\mathbf{u}_{n+1}(\mathbf{x}, t) = F(\mathbf{x}, t) - b(\mathbf{u}_n, \mathbf{u}_n)(\mathbf{x}, t) \tag{16}$$

where $F(\mathbf{x}, t) = T_t \mathbf{u}_0(\mathbf{x}) + \int_0^t T_s \mathcal{P} \mathbf{g}(\mathbf{x}, t-s) ds$, $b(\mathbf{u}, \mathbf{v}) := \int_0^t \mathcal{P} T_{t-s} \nabla \cdot (\mathbf{u} \otimes \mathbf{v}) ds$, $\mathbf{u}^{(0)}(\mathbf{x}, t) = T_t \mathbf{u}_0(\mathbf{x})$ and $\mathbf{u}_1(\mathbf{x}, t) = F(\mathbf{x}, t) + b(\mathbf{u}^{(0)}, \mathbf{u}^{(0)})(\mathbf{x}, t)$. The convergence of the iterates follows from showing that that these are iterates of a contraction on a ball in a suitable function space; c.f. Meyer (2004) and references therein.

It may be noted that the function space $L^3(\mathbf{R}^3)$, or more generally $L^n(\mathbf{R}^n)$, enjoys the very special norm scale invariance of the form

$$\|f_\lambda\|_{L^3} = \|f\|_{L^n}, \qquad f_\lambda(\mathbf{x}) = \frac{1}{\lambda} f(\lambda\mathbf{x}), \mathbf{x} \in \mathbf{R}^n, \lambda > 0. \tag{17}$$

Such Banach spaces are referred to as *limit spaces*. Excellent surveys of contemporary literature which grew out of Kato's approach are given in Meyer [42],

Cannone [10, 9]. A more general notion of limit space arises in the context of majorizing kernels and will be given in Section 3; see (64). The Picard iteration approach may also be viewed as a starting point for the probabilistic methods initiated by LeJan and Sznitman [37], though perhaps not in a way one might expect on first glance. Their's is based on a branching representation in Fourier space in which the non-linearity in the Fourier-transformed version of (13) is associated with a "binary branching product" of certain random walks in Fourier frequency space. The essential difference in the two approaches is in the nature of bounds sought for convergence of two distinct but closely related limits. In the former one seeks contraction bounds for convergence of the iteration scheme, while in the latter one seeks integrability bounds on stochastic variables, defined by a.s. limits of a stochastic iteration, to ensure finiteness of their expected values. To highlight the close connections between these two approaches it may be noted that the aforementioned uniqueness proof by Furioli, Lemarié-Rieusett and Terraneo [22] for Kato's solution in $L^3(\mathbf{R}^3)$ was obtained by arguments, specifically estimates on the associated bilinear form $b(\mathbf{u}, \mathbf{v})$, which the authors explicitly acknowledge to have been inspired by the probabilistic approach to uniqueness for a certain other space obtained by LeJan and Sznitman [37]; namely the space $\mathcal{F} = \{\mathbf{u} \in \mathcal{S}'(\mathbf{R}^3) : \sup_{\xi} |\xi|^2 |\hat{u}(\xi)| < \infty\}$ or, more accurately, the solenoidal distributions belonging to this space $\mathcal{F}$. A distribution whose Fourier transform is essentially bounded is referred to as a *pseudo-measure* in harmonic analysis. Thus an alternative description of this space is one of locally integrable tempered distributions for which $\Delta\mathbf{u}$ is a pseudo-measure; now often denoted $\mathcal{PM}^2$ where the superscript 2 refers to the $|\xi|^2$. pre-factor, c.f. Cannone [9]. In fact this may also be shown to be a Besov space; see Cannone and Planchon [11]. For the unitiated, a brief historical orientation to the various such function spaces which arise in this context may be found in Peetre [48].

Since, loosely speaking, algebraic products transform to convolution products of Fourier transforms and since derivatives transform by simple algebraic factors, the transformation of the non-linearity under Fourier transform in the (quasi-linear) Navier-Stokes equation is deceptively similar to the convolution term one obtains in Fourier space for (semi-linear ) reaction-diffusion equations, e.g. KPP. But one must not forget about the incompressibility condition, which is completely absent for the case of reaction-diffusion. Nonetheless, it is natural to search for branching Brownian motions in physical space; e.g. see McKean [40] in the case of KPP. In the case of the scalar Burgers equation, again devoid of a notion of incompressibility, one may identify such a representation; see Waymire [56], Ramirez [49]. Moreover, recently a probabilistic representation involving semi-Markovian branching Brownian motions in physical space was obtained by Ossiander [47] in which some basic notions from probabilistic potential theory make a natural appearance. In the next section we will further indicate natural probabilistic ties to harmonic analysis in the way of the "background radiation" process in physical space introduced in Gundy and Varapoulos [27] and Gundy and Silverstein [26].

The previously noted basic results capture the essence of the paradigm furnished by iterative methods in a wide range of function spaces. Namely, one either obtains for sufficiently small initial data $\mathbf{u}_0$ existence and uniqueness of global solutions, i.e. for all times $T > 0$, or one proves existence and uniqueness for arbitrarily large initial data $\mathbf{u}_0$ for short durations $T \equiv T(\mathbf{u}_0) < \infty$. A concisely formulated point for caution concerning semi-group approaches was provided in a paper by Montgomery-Smith [44] which suggests a much deeper role for incompressibility than has thus far been captured by these methods. However there are many significant questions beyond existence and uniqueness, e.g. numerical computation and Monte-Carlo simulation algorithms, boundary value problems, and steady-states, which stand to benefit by continued development of probabilistic methods. It is manifestly clear that the introduction of probabilistic reasoning has certainly lead to new insights into the structure of these equations and their solutions. On the other hand, it is also clear that incompressibility is being undervalued in these analyses to date.

## 2. Vorticity, Pressure, Incompressibility and Background Radiation

One need do no more than point to the classic papers on diffusions by W. Feller, to name just one example, to highlight and exemplify the impact that the analytic theory of linear (parabolic) partial differential equations has had on the development of probability theory. On the other hand there is no question to the significant counterpoint that the subsequent development of martingale theory, the Itô calculus and the Malliavin calculus has provided in the study of such pde's. However, unlike the fundamental equations governing diffusion, the historic origins of the Navier-Stokes equations are devoid of explicit probabilistic reasoning for their formulation - making the appearance of probabilistic methods here all the more mysterious from the outset.

There is that Laplacian term! One of the earliest applications of probabilistic ideas in connection with incompressible Navier-Stokes was introduced by Chorin [13], based on the associated *vorticity equation*. The vector cross product $\omega(\mathbf{x}, t) = \nabla \wedge \mathbf{u}(\mathbf{x}, t)$, i.e. the *curl* of the velocity field, defines the *vorticity*. Although the physical equations are derived in $n = 3$ dimensions the formulation extends to $n \geq 2$ dimensions. In $n = 1$ dimension the notion of vorticity is lost; the corresponding scalar equation for velocity $\frac{\partial u}{\partial t} = \nu \frac{\partial u}{\partial x} + u \frac{\partial u}{\partial x} + g$ being referred to as *Burgers equation*. However in the lowest dimensional case $(n = 2)$ of velocity $u = (u_1, u_2)$ in (1), one may also identify the vector vorticity $\omega \equiv (0, 0, \omega)$ with the scalar component $\omega$.

Gradients are annihilated by taking the curl since for any sufficiently smooth scalar function $a$

$$
\begin{aligned}
\nabla \wedge \nabla a &= \nabla \wedge (\partial a/\partial x_1, \partial a/\partial x_2, \partial a/\partial x_3) \\
&= (\partial^2 a/\partial x_2 \partial x_3 - \partial^2 a/\partial x_3 \partial x_2, \partial^2 a/\partial x_3 \partial x_1 - \partial^2 a/\partial x_1 \partial x_3, \\
&\quad \partial^2 a/\partial x_1 \partial x_2 - \partial^2 a/\partial x_2 \partial x_1). \\
&= 0. \tag{18}
\end{aligned}
$$

In particular taking the curl provides another approach to remove the (unknown) pressure gradient in (1). Applying some simple vector identities, e.g. see Chorin and Marsden [15, pp. 163-164], namely $(\mathbf{u} \cdot \nabla)\mathbf{u} = \frac{1}{2}\nabla(|\mathbf{u}|^2) + \omega \wedge \mathbf{u}$ and $\nabla \wedge (\omega \wedge \mathbf{u}) = \omega(\nabla \cdot \mathbf{u}) - \mathbf{u}(\nabla \cdot \omega) + (\mathbf{u} \cdot \nabla)\omega - (\omega \cdot \nabla)\mathbf{u}$, and using incompressibility and (18), one may express the curl of $\mathbf{u}$ satisfying (1) as

$$\frac{\partial \omega}{\partial t} = \nu\Delta\omega - (\mathbf{u} \cdot \nabla)\omega + (\omega \cdot \nabla)\mathbf{u} + \nabla \wedge \mathbf{g}, \quad \nabla \cdot \omega = 0. \tag{19}$$

A solution $\omega$ of (19) may then be used to obtain a solution $\mathbf{u}$ to (1) and, in turn, the pressure $p$ as follows: According to the *Helmholtz-Hodge decomposition* any smooth vector field $\mathbf{u}$ which decays sufficiently fast at infinity may be uniquely represented as the superposition of a gradient and a (necessarily incompressible) curl

$$\mathbf{u} = \mathbf{u}_1 + \mathbf{u}_2, \tag{20}$$

where

$$\mathbf{u}_1 = -\nabla\Phi, \quad \mathbf{u}_2 = \nabla \wedge \Psi \tag{21}$$

for a scalar potential $\Phi$ and vector potential $\Psi$ obtained by solving a (vector) Poisson equation

$$\Delta\Phi = -\nabla \cdot \mathbf{u}, \qquad \Delta\Psi = -\Delta \wedge \mathbf{u}, \quad \nabla \cdot \Psi = 0. \tag{22}$$

In particular one may recover $\mathbf{u}$ from vorticity $\omega$ via the following so-called *Biot-Savart law*:

$$\mathbf{u}(\mathbf{x}, t) = L * \omega(\mathbf{x}, t), \tag{23}$$

where $L$ is a linear (convolution) kernel furnished by the divergence free solution $\Psi$ to the Poisson equation, i.e. by inverting the Laplacian via (for $n = 3$)

$$\Psi(\mathbf{x}, t) = \Delta^{-1}(-\omega) = -\frac{1}{4\pi}\int_{\mathbf{R}^3}\frac{1}{|\mathbf{x} - \mathbf{y}|}\omega(\mathbf{y}, t)d\mathbf{y}, \quad \nabla \cdot \Psi = 0, \tag{24}$$

and hence (for $n = 3$)

$$\begin{aligned}\mathbf{u}(\mathbf{x}, t) &= \nabla \wedge \left(-\frac{1}{4\pi}\int_{\mathbf{R}^3}\frac{1}{|\mathbf{x} - \mathbf{y}|}\omega(t, \mathbf{y})d\mathbf{y}\right) = \frac{1}{4\pi}\int_{\mathbf{R}^3}\nabla_{\mathbf{x}}\frac{1}{|\mathbf{x} - \mathbf{y}|}\omega(t, \mathbf{y})d\mathbf{y}\\ &= L * \omega(\mathbf{x}, t).\end{aligned} \tag{25}$$

Similarly for $n \neq 3$ one needs only to apply the appropriate potential kernel to invert the Laplacian. In particular for $n = 2$ dimensions one has

$$\mathbf{u}(\mathbf{x}, t) = \nabla \wedge \left(-\frac{1}{2\pi}\int_{\mathbf{R}^2}\log(|\mathbf{x} - \mathbf{y}|)\omega(\mathbf{y}, t)d\mathbf{y}\right) = L * \omega(t, \mathbf{x}), \tag{26}$$

where

$$L(\mathbf{x}) = \frac{1}{2\pi|\mathbf{x}|^2}(x_2, -x_1), \quad \mathbf{x} = (x_1, x_2) \in \mathbf{R}^2. \tag{27}$$

For numerical computations Chorin [13] viewed (19) by an *operator splitting*, i.e. along the lines of a Trotter-Kato semi-group product, into successive incremental advective and diffusive flow processes using distinct numerical schemes.

In this connection, assuming $n = 2$ observe that the term $\omega \cdot \nabla$ in (19) vanishes for $n = 2$; here one should recall the convention of using the scalar $\omega$ and vector $(0, 0, \omega)$ interchangeably in $n = 2$ dimensions. Thus, assuming also $g \equiv 0$, the two-dimensional vorticity equation describes the evolution of a scalar $\omega = \omega(\mathbf{x}, t)$ by

$$\frac{\partial \omega}{\partial t} = \nu \Delta \omega - (\mathbf{u} \cdot \nabla)\omega. \tag{28}$$

Now, if $\mathbf{u}$ were given (computed) over an increment of time, then the scalar component of vorticity would appear to be governed by a diffusion equation. In this regard, denote two-dimensional standard Brownian motion by $\mathbf{W} = \{(W_1(t), W_2(t)) : t \geq 0\}$ and let $\mathbf{Z} := \{(Z_1(t), Z_2(t)) : t \geq 0\}$ be the diffusion defined by the associated stochastic differential equation

$$d\mathbf{Z}(t) = \sqrt{2\nu} d\mathbf{W}(t) - \mathbf{u}(\mathbf{Z}(t), t)dt \tag{29}$$

Then one is naturally led to consider

$$\omega(\mathbf{x}, t) := E\omega_0(\mathbf{Z}(t)). \tag{30}$$

This is the starting point for an approach which has been exploited heavily as a numerical Monte-Carlo method, and has been studied closely by analysts and probabilists alike; e.g. Chorin [13], Marchioro and Pulverenti [39] Goodman [24], Long [38], Szumbarski and Wald [51], Meleard [41]. In this regard, Busnello, Flandoli, and Romito [8] also obtain a natural probabilistic version of the Biot-Savart law in three-dimensions. While we will not pursue these particular approaches any further, we will next take the idea of considering linearization and Duhamels principle in an apparently new direction which exploits incompressibility in an interesting way.

To recover the pressure from its projection onto a divergence free vector field let us recall (12). The matrix elements $\frac{\xi_i}{|\xi|} \frac{\xi_j}{|\xi|}$ may be be viewed as Fourier symbols of iterated Riesz transforms $R_i R_j$, i.e. $\widehat{R_j f}(\xi) = -\frac{\xi_j}{|\xi|} \hat{f}(\xi)$; note the signs and constants depend on the version (10) of Fourier transform used. In particular, on appropriate function spaces, e.g. $L^2(\mathbf{R}^3)^3$, one has

$$\mathcal{P} = \mathbf{I} + \mathcal{R}, \tag{31}$$

where $\mathcal{R} = ((R_i R_j))_{1 \leq i, j \leq 3}$ is the matrix of iterated Riesz transforms. To exploit incompressibility in (1) note that the terms linear in $\mathbf{u}$ are divergence free. Thus take the divergence to annihilate the linear terms $\frac{\partial \nabla \cdot \mathbf{u}}{\partial t} = 0 = \Delta(\nabla \cdot \mathbf{u}) = \nabla \cdot (\Delta \mathbf{u})$, and assume no forcing $\mathbf{g} = 0$, so that one has

$$\nabla \cdot (\mathbf{u} \cdot \nabla)\mathbf{u} = -\nabla \cdot \nabla p = -\Delta p. \tag{32}$$

Taking Fourier transforms with continued use of incompressibility one arrives at the following determination of the pressure

$$p = \sum_{1 \leq j, k \leq 3} R_j R_k(u_j u_k). \tag{33}$$

Apart from what may be gleaned from its potential theoretic role in making projections onto divergence free vector fields, a basic understanding of incompressibility remains illusive from a probabilistic perspective. However, even in the case of *linearized Navier-Stokes* equations, the incompressibility condition plays a fascinating role. Two common linearizations are (i.)*Stokes* equations and (ii.)*Oseen* equations, respectively given by linearization around zero velocity or by linearization around a non-zero constant $\overline{\mathbf{U}}$, respectively. Specifically one has, respectively, for Stokes

$$\frac{\partial \mathbf{u}}{\partial t} = \nu \Delta \mathbf{u} - \nabla p + \mathbf{g}, \quad \nabla \cdot \mathbf{u} = 0. \tag{34}$$

and for Oseen

$$\frac{\partial \mathbf{u}}{\partial t} + \overline{\mathbf{U}} \cdot \nabla \mathbf{u} = \nu \Delta \mathbf{u} - \nabla p + \mathbf{g}, \quad \nabla \cdot \mathbf{u} = 0. \tag{35}$$

For simplicity and without any significant loss of generality let us consider the Stokes equation (34). Observe that for given twice continuously differentiable and *incompressible* initial data $\mathbf{u}_0 = (u_0^{(1)}, u_0^{(2)}, u_0^{(3)})$, one has that in the unforced case ($\mathbf{g} \equiv 0$) the solution $(\mathbf{u}, p) = ((u_1, u_2, u_3), p)$ is given by the familiar semi-group representation

$$u_j(\mathbf{x}, t) = E\mathbf{u}_0^{(j)}(\sqrt{2\nu}\mathbf{W}(t)), \; j = 1, 2, 3, \quad p \equiv c, \tag{36}$$

where $\mathbf{W} = \{(W_1(t), W_2(t), W_3(t)) : t \geq 0\}$ is three-dimensional standard Brownian motion, and $c$ is an arbitrary constant. In particular the linearity in $\mathbf{u}$ and the *apparent affine* form of (34) due to $-\nabla p$ are reconciled by the incompressibility. In other words, the incompressibility of $\mathbf{u}$ (at times $t > 0$) is directly inherited from that of $\mathbf{u}_0$, causing $-\nabla p$ to instantaneously "drop out" at $t = 0$. Accordingly one may view the matrix-valued Gaussian (heat) kernel

$$\mathbf{K}(\mathbf{z}, t) = \mathbf{diag}[k(\mathbf{z}, t), k(\mathbf{z}, t), k(\mathbf{z}, t)],$$

$$k(\mathbf{z}, t) = (4\pi\nu t)^{-\frac{3}{2}} e^{-\frac{|\mathbf{z}|^2}{4\pi\nu t}} \tag{37}$$

furnished by the transition probabilities of the Brownian motion, as the fundamental solution to (34) for incompressible initial data.

Let us consider, on the other hand, the evolution of smooth but compressible initial data $\mathbf{u}_0$. Also let us include a possibly compressible forcing term $\mathbf{g}$. We seek an incompressible vector field $\mathbf{u}(\mathbf{x}, t) \equiv \mathcal{P}\mathbf{u}(\mathbf{x}, t), t > 0$, such that by orthogonality (6) of $\nabla p$ to incompressible $\mathbf{u}$ for $t > 0$,

$$\frac{\partial \mathcal{P}\mathbf{u}}{\partial t} = \nu \Delta \mathcal{P}\mathbf{u} + \mathcal{P}(\mathbf{g}), \; t > 0, \quad \mathcal{P}\mathbf{u}(\mathbf{x}, 0^+) = \mathcal{P}\mathbf{u}_0(\mathbf{x}). \tag{38}$$

The compressible initial data $\mathbf{u}_0$ will necessitate that $\mathbf{u}(\mathbf{x}, t) \equiv \mathcal{P}\mathbf{u}(\mathbf{x}, t)$ for $t > 0$, be discontinuous at $t = 0$. Applying Duhamel's principle to (38) one

arrives at

$$
\begin{aligned}
\mathbf{u}(\mathbf{x},t) &= \int_{\mathbf{R}^3} \mathbf{K}(\mathbf{z}-\mathbf{x},t)\mathcal{P}\mathbf{u}_0(\mathbf{z})d\mathbf{z} + \int_0^t \int_{\mathbf{R}^3} \mathbf{K}(\mathbf{x}-\mathbf{z},t-s)\mathcal{P}\mathbf{g}(\mathbf{z},s)ds d\mathbf{z} \\
&= \int_{\mathbf{R}^3} \mathbf{\Gamma}(\mathbf{z}-\mathbf{x},t)\mathbf{u}_0(\mathbf{z})d\mathbf{z} + \int_0^t \int_{\mathbf{R}^3} \mathbf{\Gamma}(\mathbf{x}-\mathbf{z},t-s)\mathbf{g}(\mathbf{z},s)ds d\mathbf{z}, \quad (39)
\end{aligned}
$$

where

$$
\mathbf{\Gamma} = \{\mathbf{I} + \mathcal{R}\}(\mathbf{K}), \tag{40}
$$

and $\mathcal{R} = ((R_j R_k))_{1 \le j,k \le 3}$ is the matrix of iterated Riesz transforms. The simplest way to check the commutativity required for the second equation in the display (39) is by Fourier transforms along the lines indicated earlier for (31). To compute the pressure one simply takes the divergence in (34) and use incompressibility to get

$$
\nabla \cdot \nabla p(\mathbf{x},t) = \Delta p = \nabla \cdot \mathbf{g}, \ t > 0. \tag{41}
$$

For $t = 0$ one may apply the Helmholtz-Hodge decomposition to recover $u_0$ from $\mathcal{P}\mathbf{u}_0$ defining the pressure at $t = 0$ by

$$
p(\mathbf{x},0) = \Phi = \Delta^{-1}(-\nabla \cdot \mathbf{u}_0). \tag{42}
$$

In particular, therefore,

$$
p(\mathbf{x},t) = \Delta^{-1}(-\nabla \cdot \mathbf{u}_0)(\mathbf{x})\delta_0(t) + \Delta^{-1}(\nabla \cdot \mathbf{g})(\mathbf{x},t) \tag{43}
$$

One may note that (39) actually embodies two Duhamel principles; one applied and one derived! Namely, we apply a Duhamel principle for Brownian motion to (38) with $\mathcal{P}(\mathbf{g})$ as the forcing, and derive a Duhamel principle for (34) with arbitrary initial data using the fundamental solution (40) of the unforced linear incompressible pressurized equation and adding $\mathbf{g}$ as the forcing.

The representation (40) of the (signed) fundamental solution $\mathbf{\Gamma}$ of Stokes equation for $(\mathbf{u}, p)$ is the well-known classic formula of Oseen [46]; also see Solonnikov [50], Thomann and Guenther [53]. One may check a semi-group property of the form

$$
\mathbf{\Gamma}(\mathbf{x},t+s) = \int_{\mathbf{R}^3} \mathbf{\Gamma}(\mathbf{x}-\mathbf{y},t)\mathbf{\Gamma}(\mathbf{y},s)d\mathbf{y}. \tag{44}
$$

Moreover, $\mathbf{\Gamma}$ has the following equivalent representation

$$
\mathbf{\Gamma} = \mathbf{K} + \mathbf{Hess}(\psi), \quad \Delta\psi = -\mathbf{K}, \tag{45}
$$

where $\mathbf{Hess}$ denotes the *Hessian matrix*. A simple way to check this is via Fourier transforms since $\frac{\xi_j \xi_k}{|\xi|^2}$ is the Fourier symbol for the operator $\frac{\partial^2}{\partial x_j x_k}\Delta^{-1}$. Thomann and Guenther [53] recently made an explicit computation of the full fundamental solution in terms of special functions which has had major implications for analysis of Stokes equation in unrestricted domains. Their precise formulae may be expressed as

$$\Delta^{-1} K = \frac{1}{4\pi|\mathbf{z}-\mathbf{x}|} Erf(\frac{|\mathbf{z}-\mathbf{x}|}{\sqrt{4t}})$$

$$-R_j R_k(K) = K(\mathbf{z}-\mathbf{x}, t)\{-\frac{1}{3} \, {}_1F_1(1, \frac{5}{2}, \frac{|\mathbf{z}-\mathbf{x}|^2}{4t})\delta_{jk}$$

$$+\frac{(x_j - z_j)(x_k - z_k)}{|\mathbf{z}-\mathbf{x}|^2}\left({}_1F_1(1, \frac{5}{2}, \frac{|\mathbf{z}-\mathbf{x}|^2}{4t}) - 1\right)\}, \qquad (46)$$

where $\mathrm{Erf}(s) = \frac{2}{\sqrt{\pi}}\int_0^s e^{-y^2}dy$ denotes the error function and ${}_1F_1(a, c, u) = \sum_{n=0}^{\infty} \frac{(a)_n}{(c)_n n!} u^n$, for $(b)_n = \frac{\Gamma(b+n)}{\Gamma(b)}$, is a confluent hypergeometric function. This is a hugely significant formula in view of the rich mathematical literature on parametric dependencies, asymptotic growth and/or decay properties etc of the special functions which appear in the formula (46).

The matrix $\mathbf{\Gamma}(\mathbf{x}, 0)$ is the orthogonal projection of the (vector) Dirac distribution $\delta_{\mathbf{x}}$ onto divergence-free vector fields. Thus, the effect of incompressibility is instantaneous in that the solution is obtained by projecting the initial data onto divergence-free vector fields and then solving (34) with incompressible initial data as above. In fact one may observe that

$$\mathbf{\Gamma}(\mathbf{x}, t) = \int_{\mathbf{R}^3} \mathbf{K}(\mathbf{x}-\mathbf{y}, t)\mathbf{\Gamma}(\mathbf{y}, 0)d\mathbf{y}, \quad \mathbf{x} = (x, y, z) \in \mathbf{R}^3. \qquad (47)$$

The transformation $\mathbf{\Gamma}(\mathbf{x}_0, 0)\mathbf{u}_0$ is the (instantaneous) projection of $\mathbf{u}_0$ onto its divergence free (incompressible) component.

The sign changes and discontinuity at $t = 0$ in the fundamental solution make the issue perplexing from a probabilistic standpoint. However, there is some hope that a probabilistic unraveling of this semi-group might lead to ways to deal with boundaries. Thomann and Guenther (personal communication) have pointed out that if one poses a Neumann boundary condition on the half-space, the usual "reflection method" applied to the fundamental solution $\Gamma$ for the unrestricted problem fails.

In view of the prominent role of Riesz transforms it seems natural to explore a role for Gundy-Varapoulos-Silverstein's *background radiation* process in the present framework. Thus we close this section with elements of such a new connection. The background radiation process $\mathbf{Z} = \{Z_t \equiv (Y_t, \mathbf{X}_t) : -\infty < t \leq 0\}$ belongs to the class of "approximate Markov processes" defined, as a collection of measurable functions, on a sigma-finite measure spaces $(\Omega, \mathcal{F}, P)$ taking values in the half-space $H_3 = [0, \infty) \times \mathbf{R}^3$. Foregoing its construction, the process is defined by the following properties (Gundy [25]) :

1  $\{Y_t : -\infty < t \leq 0\}$ is distributed as a one-dimensional Brownian motion on $(0, \infty)$ started at (arbitrary) $y_0 > 0$ and killed upon entry to the state 0.

2  $\lim_{t \uparrow 0} Y_t = 0$ a.s.

3 $\{\mathbf{X}_t : -\infty < t \leq 0\}$ is an independent Brownian motion with initial distribution (at time "$t = -\infty$") given by Lebesgue measure on $\mathbf{R}^3$ in the sense that for each $y > 0$, defining $\tau_y = \inf\{t \leq 0 : Y_t = y\}$, the process $\{(y, \mathbf{X}_{\tau_y + t}) : 0 \leq t \leq -\tau_y\}$ is distributed as a Brownian motion started with Lebesgue measure on the hyper-plane $\{y\} \times \mathbf{R}^3$ at $t = 0$.

4 $\mathbf{X}_0 = \lim_{t \uparrow 0} \mathbf{X}_t$ exists a.s. and is distributed as Lebesgue measure on $\mathbf{R}^3$, i.e. the induced measure $P \circ \mathbf{X}_0^{-1}$ is Lebsesgue measure.

These properties were used by Gundy and Varapoulos [27] in Gundy and Silverstein [26] to prove the following formula for iterated Riesz transforms of a function $f$ belonging to Schwartz space $\mathcal{S}$

$$R_j R_k f(x) = bE(\int_{-\infty}^{0} A_{jk} \nabla f_0 \cdot dZ_t | Z_0 = (0, x)) \tag{48}$$

where $f_0$ denotes the harmonic extension of the function $f$ on the $\mathbf{R}^3$ to the half-space $H_3$, the matrix $A_{jk} = \mathbf{e}_j \otimes \mathbf{e}_k$, i.e. having 1 in location $(j, k)$ and 0's elsewhere, defines the indicated martingale transform. With a little more work one may in fact calculate $b = -\frac{1}{2}$ using their indicated methods. In view of (40) and viewing conditional expectations as disintegration formulae, one may express the solution to Stokes linearization (34) in terms of background radiation as follows, writing $\mathbf{u} = (u^{(1)}, u^{(2)}, u^{(3)})$:

$$
\begin{aligned}
u^{(j)}(\mathbf{x}, t) &= \int_{\mathbf{R}^3} k(\mathbf{z} - \mathbf{x}, t) \mathbf{u}_0^{(j)}(\mathbf{z}) d\mathbf{z} \\
&+ \sum_{k=1}^{3} R_j R_k(K)(\mathbf{z} - \mathbf{x}, t) \mathbf{u}^{(k)}(\mathbf{z}) d\mathbf{z} \\
&= E_{X_0 = \mathbf{x}} u_0^{(j)}(\mathbf{X}_t^*) \\
&+ \sum_{k=1}^{3} \int_{\mathbf{R}^3} E(\theta^{(jk)}(\mathbf{x}, t) | \mathbf{Z}_0 = (0, \mathbf{z})) u_0^{(k)}(\mathbf{z}) d\mathbf{z},
\end{aligned}
\tag{49}
$$

where $\{X_t^* : t \geq 0\}$ is the *reversed-time* process (Brownian motion) introduced in Gundy and Silverstein [26], and $\theta^{(jk)}$ is the martingale transform defined by

$$\theta^{(jk)} \equiv \theta^{(jk)}(\mathbf{x}, t) = -\frac{1}{2} \int_{-\infty}^{0} \mathbf{e}_j \otimes \mathbf{e}_k \nabla K_0(\mathbf{Z}_s) \cdot d\mathbf{Z}_s \tag{50}$$

and for fixed $t, \mathbf{x}$, $K_0$ is the harmonic extension of the the heat kernel $\mathbf{z} \rightarrow K(\mathbf{z} - \mathbf{x}, t)$ to the half-space $H_3$. In particular, defining the matrix-valued function

$$\Theta(\mathbf{x}, t) = ((\theta^{(jk)}(\mathbf{x}, t)))_{1 \leq j, k \leq 3}, \tag{51}$$

one has the representation

$$\mathbf{u}(\mathbf{x}, t) = E_{\mathbf{Z}_0 = (0, \mathbf{x})} \mathbf{u}_0(\mathbf{X}_t^*) + E(\Theta(\mathbf{x}, t) \mathbf{u}_0(\mathbf{Z}_0)). \tag{52}$$

Form here one may regard the non-linear term in (1) as a forcing and exploit a Duhamel principle to obtain an integral equation involving background radiation. While this connection has not been pursued beyond this survey article, an alternative real-space semi-Markov cascade representation of Navier-Stokes is fully developed in Ossiander [47] which also contains intriguing connections with the fundamental solution obtained by Thomann and Guenther [53] for linearized flow.

## 3. Multiplicative Stochastic Cascade

One may give purely fluid mechanical derivations of diffusion equations based on mass conservation and linear flux laws such as Ficks and Fourier laws, however the most compelling physical arguments are no doubt Einstein's implicit use of the central limit theorem and law of large numbers. While the modern approach of hydrodynamic scaling limits of interacting particle systems seeks to arrive at the Navier-Stokes equations from a more intrinsic probabilistic or statistical mechanical foundation , such approaches are fundamentally much more complex than the connections that are possible between Brownian motion and the heat equation; see e.g. Esposito, Marra and Yau [16], Yau [55]. That said, owing to a remarkable insight of LeJan and Sznitman [37] it is nevertheless possible to approach certain classes of solutions to these equations by probabilistic constructions of branching random walks and multiplicative cascades in a way which preserves some basic Markovian semi-group structure. This is the topic to be briefly described in this section.

As remarked earlier the multiplicative cascade approach appears as a mathematical device, closely related to Kato's contraction map techniques for obtaining a mild solution on appropriate function spaces. While a pure analyst may take solace in this realization, as a probabilistic technique it has additional virtues. Chief among these is the "natural" emergence of solution spaces (and their attendant norms) for which one may obtain unique global solutions under a "small" initial data constraint as a consequence of the probabilistic representation. The calculation of steady state solutions are also quick to follow in these function spaces; see Bhattacharya et al. [7]. New questions and ideas with regard to Monte Carlo methods emerge naturally in more pragmatic considerations of numerical solutions; see Ramirez [49]. However, whether this approach can lead to truly significant new results on the existence and uniqueness of global solutions under smooth initial data without regard to "size" remains to be seen. It is clear that a more substantial use of incompressibility is lacking in the existing theory; c.f. Montgomery-Smith [44].

**Example** (Warm-up). Some of the basic ideas under multiplicative cascade representations are more generally applicable to diverse classes of evolution equations, including certain linear parabolic and fractional diffusion equations, semi-linear reaction-diffusions and some quasi-linear equations such as two-dimensional incompressible Navier-Stokes equations and one-dimensional Burgers' equation; see Bhattacharya et al. [7]. Although interesting from a number of other points of view, e.g. physical space representations, Monte-Carlo simulations, etc., these examples are not capable of furthering understanding of incompressibility. That said,, the following extremely simple example serves to illustrate some of the most basic graph theoretic and probabilistic ideas involved in this approach. Consider

$$u_t = a\Delta u + b \cdot \nabla u \quad u(x,0) = u_0(x), \tag{53}$$

in $n \geq 1$ dimensions, where $a > 0$, and $b \in \mathbf{R}^n$ are constants. To quickly get the flavor of the method, consider the spatial Fourier transform of (53)

$$\hat{u}(\xi,t) = \hat{u}_0(\xi)e^{-a|\xi|^2 t} + \frac{ib \cdot \xi}{a|\xi|^2} \int_0^t a|\xi|^2 e^{-a|\xi|^2 s} \hat{u}(\xi,t-s)ds. \tag{54}$$

Now consider the random linear tree $\tau_\theta(t)$ rooted at a vertex $\theta$ of type $\xi_\theta = \xi$ which, after an exponential length of time is replaced by a single vertex $\langle 1 \rangle$ of the same type $\xi_{\langle 1 \rangle} = \xi$. Proceeding in this manner one may calculate that the very familiar solution $\exp(-a|\xi|^2 + ib \cdot \xi)\hat{u}_0(\xi)$ is the expectation of the random product $\curlyvee(\theta,t)$ initialized by $\xi_\theta = \xi$ and consisting of factors $m(\xi) = ib \cdot \xi/a|\xi|^2$ at each vertex until termination where one attaches the end factor $\hat{u}_0(\xi)$. The indicated stochastic recursion then yields

$$\curlyvee(\theta,t) = m(\xi)^{N(t)}\hat{u}_0(\xi),$$

and the solution is represented as

$$\hat{u}(\xi,t) = E_{\xi_\theta=\xi}\curlyvee(\theta,t) = E_{\xi_\theta=\xi}m(\xi)^{N(t)}\hat{u}_0(\xi), \tag{55}$$

where $N(t)$ is the Poisson process with parameter $\lambda(\xi) := a|\xi|^2$ which counts the number of times the exponential clocks ring before time $t$. In particular the Poisson process formally suggests a Fourier *dual* role to that played by the standard Brownian motion in the real space expectation formula. Similarly one may obtain a formal "dual" Feynman-Kac formula (in the sense of this Fourier transform representation) under the *complex measure condition* on coefficients given by Itô [29]; see Kolokoltsov [31], and Chen et al. [12]. In particular this approach makes Itô's complex measure condition completely natural from a probabilistic point of view. One may also obtain a Fourier version of McKean's [40] branching Brownian motion formula for KPP, and other interesting equations; see Orum [45]. Also included are, for example, the generalized fractional Burgers equation of the type considered by Woyczynski, Biler, and Funaki [54].

The essential idea introduced by LeJan and Sznitman [37] and presented here in the generalized form developed in Bhattacharya et al. [6, 7] is as follows. Without loss of generality assume that the forcing $\mathbf{g}$ is incompressible; else replace $\mathbf{g}$ by $\mathcal{P}\mathbf{g}$ below. Either by directly taking Fourier transforms in the projected equation (13) and then introducing an exponential integrating factor, or by taking Fourier transforms in (14), one may deduce the following integral equation in Fourier frequency space:

$$\hat{\mathbf{u}}(\xi, t) = e^{-\nu|\xi|^2 t}\hat{\mathbf{u}}_0(\xi) + \int_0^t e^{-\nu|\xi|^2 s}\big\{|\xi|(2\pi)^{-\frac{3}{2}}$$

$$\int_{\mathbf{R}^3} \hat{\mathbf{u}}(\eta, t-s) \otimes_\xi \hat{\mathbf{u}}(\xi - \eta, t-s)d\eta + \hat{\mathbf{g}}(\xi, t-s)\big\}ds, \quad \text{(FNS)}$$

where for complex vectors $w, z$

$$w \otimes_\xi z = -i(e_\xi \cdot z)\pi_{\xi^\perp} w, \ e_\xi = \frac{\xi}{|\xi|}, \ \text{and } \pi_{\xi^\perp} w = w - (e_\xi \cdot w)e_\xi \qquad (55)$$

is the projection of $w$ onto the plane orthogonal to $\xi$, and $\nu > 0$ is the viscosity parameter. For $\xi \neq 0$, LeJan and Sznitman [37] rescale the equation (FNS) to normalize the integrating factor $e^{-\nu|\xi|^2 s}$ to the exponential probability density $\nu|\xi|^2 e^{-\nu|\xi|^2 s}$. They then observe that the resulting equation is precisely the form for a branching random walk recursion for $\chi(\xi, t) := \nu|\xi|^2\hat{\mathbf{u}}(\xi, t)$, for which the kernel $|\xi - \eta|^{-2}|\eta|^{-2}$ is naturally constrained by integrability to dimensions $n \geq 3$ for normalization to a transition probability. Alternatively, Bhattacharya et al. [6] introduce a Fourier multiplier $1/h$, where $h(\xi), \xi \in W_h := \{\xi \in \mathbf{R}^3 : \xi \neq 0\}$ is a positive function such that $0 < h * h(\xi) < \infty$, which is used to re-scale the Fourier transformed equation (FNS) by factors $1/h(\xi)$ in compensation for the *temporal* normalization of the exponential factor to a probability. Namely, we consider the equation (FNS$_h$) obtained from (FNS) as

$$\chi(\xi, t) = e^{-\nu t|\xi|^2}\chi_0(\xi) + \int_0^t \nu|\xi|^2 e^{-\nu t|\xi|^2 s}\big\{\frac{1}{2}m(\xi)\int_{W_h \times W_h} \chi(\eta_1, t-s)$$

$$\otimes_\xi \ \chi(\eta_2, t-s)H(\xi, d\eta_1 \times d\eta_2) + \frac{1}{2}\varphi(\xi, t-s)\big\}ds, \quad \xi \in W_h \ \ \text{(FNS}_h\text{)}$$

where

$$m(\xi) = \frac{2h * h(\xi)}{\nu(2\pi)^{\frac{3}{2}}|\xi|h(\xi)}, \quad \chi_0(\xi) = \frac{\hat{u}_0(\xi)}{h(\xi)}, \quad \varphi(\xi, t) = \frac{2\hat{g}(\xi, t)}{\nu|\xi|^2 h(\xi)}, \qquad (56)$$

and $H(\xi, d\eta_1 \times d\eta_2)$ is for $\xi \in W_h$ the transition probability kernel, with support contained in the set $\{(\eta_1, \eta_2) \in W_h \times W_h : \eta_1 + \eta_2 = \xi\}$, defined by

$$\int_{W_h \times W_h} f(\eta_1, \eta_2)H(\xi, d\eta_1 \times d\eta_2) = \int_{W_h} f(\xi - \eta, \eta)\frac{h(\xi - \eta)h(\eta)}{h * h(\xi)}d\eta \qquad (57)$$

for bounded, Borel measurable $f : W_h \times W_h \to \mathbf{R}$. In view of the *conservation of wave-numbers* $\eta_1 + \eta_2 = \xi$ reflected in the support of the transition probability kernel, the definition of the supporting set of frequencies $W_h$ may be relaxed to a semi-group $W_h \subseteq \{\xi \in \mathbf{R}^3 : \xi \neq 0\}$. However in general we include the following additional *exterior condition* in defining (FNS$_h$):

$$\chi(\xi, t) = 0, \quad \xi \in W_h^c, t \geq 0. \qquad (58)$$

One may also note that with the exception of $\partial\hat{\mathbf{u}}/\partial t$ and $\hat{g}$, the Fourier transform of each term in (1) vanishes at zero frequency. Thus the value of $\hat{\mathbf{u}}(0,t) \equiv \int_{\mathbf{R}^3} \mathbf{u}(\mathbf{x},t)d\mathbf{x}$ is explicitly determined from the initial data and forcing by

$$\frac{d}{dt}\hat{\mathbf{u}}(0,t) = \hat{\mathbf{g}}(0,t), \quad \hat{\mathbf{u}}(0,0) = \hat{\mathbf{u}}_0(0). \tag{59}$$

The probabilistic approach is based upon an interpretation of the re-scaled integral equation. This is achieved in terms of expectation values of multiplicative cascade solutions to stochastic recursions generated by a multi-type branching random walks in Fourier space. The process may be initiated at the root vertex $\theta$ of a binary tree graph having wave-number type $\xi \neq 0$. After an exponential holding time $S_\theta$ with rate $\nu|\xi|^2$ an independent coin flip is made. If the outcome is tail, denoted $[\kappa_\theta = 0]$, then the tree terminates as a single edge of length $S_\theta$ connecting $\theta$ and $< 1 >$. If the outcome is head, denoted $[\kappa_\theta = 1]$, then the edge of length $S_\theta$ branches into two vertices $< 1 >, < 2 >$ of respective types $\xi_1, \xi_2, \xi_1 + \xi_2 = \xi$, selected according to the probability kernel $H(\xi, \cdot)$. Given $\xi$ two independent exponential clocks with respective rates $\nu|\xi_1|^2$ and $\nu|\xi_2|^2$ are reset at the vertices $< 1 >, < 2 >$ and the process is repeated with each of these vertices as "roots". The root $\theta$ has an associated birth time $B_\theta \equiv 0$ and an end-time $S_\theta$ at which it is terminated or it branches. Similarly each descendent vertex $v =< v_1, \ldots, v_k >, v_j \in \{1, 2\}$, has a corresponding birth time $B_v := \sum_{j=1}^{k-1} S_{<v_1,\ldots,v_j>}$ and an end-time $B_v + S_v$. At a given time $t > 0$ one assigns multiplicative weights to vertices $v$ such that either $[B_v + S_v < t, \kappa_v = 0]$ or $[B_v < t \leq B_v + S_v]$. In the first case the re-scaled forcing $\varphi(\xi_v)$ is assigned, while in the latter case the re-scaled initial data $\chi(\xi_v)$ is assigned, where $\xi_v$ denotes the frequency type of $v$; see Figure 1.

One may now use the stochastic branching model represented by the collection of random variables $\{\xi_\mathbf{v}, \kappa_\mathbf{v}, S_\mathbf{v} : \mathbf{v} \in \mathcal{V} := \cup_{k=0}^\infty \{1, 2\}^k\}$ to recursively define a random functional related to (FNS) through its expected value. Namely, for measurable functions $\chi_0 : W_h \to \mathbf{C}^3$ and $\varphi : W_h \times [0, \infty) \to \mathbf{C}^3$, and for $\xi_\theta = \xi \in W_h, t \geq 0$, the stochastic functional $\chi(\theta, t)$ is recursively defined starting from $v = \theta$ by

$$\chi(v, s) = \begin{cases} \chi_0(\xi_v), & \text{if } S_v > s \\ \varphi(\xi_v, s - S_v) & \text{if } S_v \leq s, \kappa_v = 0 \\ m(\xi_v)\chi(< v1 >, \ s - S_v) \otimes_{\xi_v} \chi(< v2 >, s - S_v), & else \end{cases} \tag{60}$$

where the product $\otimes_\xi$ and factors $m(\xi)$ are defined in (55), (56), respectively.

Since the underlying discrete parameter branching process defined by the equally likely Bernoulli offspring distribution is a critical binary Galton-Watson process, the recursion will terminate in a finite number of iterations with probability one. In particular, $\chi(\theta, t)$ is a.s. a *finite* product of values of $\chi_0$ and/or $\varphi$ evaluated at randomly selected Fourier frequencies selected by the branching

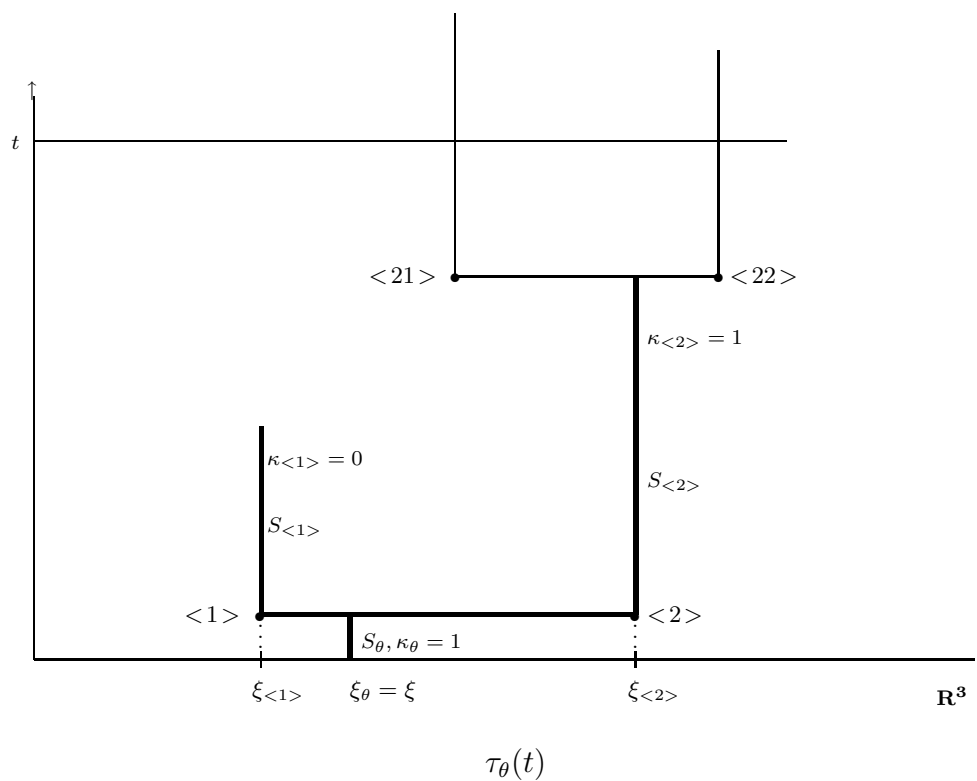

$$\tau_\theta(t)$$

Fig 1. *Schematic of a tree indexed branching random walk.*

random walk. For example, the functional evaluation of the sample tree in the figure is given by

$$\chi(\theta, t) = m(\xi_\theta) m(\xi_{<2>}) \varphi(\xi_{<1>}, t - S_{<1>} - S_\theta) \otimes_{\xi_\theta} [\chi_0(\xi_{<21>}) \otimes_{\xi_{<2>}} \chi_0(\xi_{<22>})]$$

In general $\chi(\cdot)$ is a random multiplicative functional of scalar values $m(\cdot)$ and Fourier transformed initial data and/or forcing (vector) values over the vertices of a multi-type branching random walk tree $\tau_\theta(t)$ initiated from a single progenitor of type $\xi_\theta = \xi$. In general the scalar and vector valued factors are evaluated at the wave-number (type) of the respective vertices appearing in the tree $\tau_\theta(t)$, with the initial and forcing terms appearing at the end-nodes.

This generalizes branching random walks in the sense of LeJan and Sznitman [37] for $n = 3$ dimensions where $h(\xi) = |\xi|^{-2}, \xi \in W_h = \mathbf{R}^3 \backslash \{0\}$. The essential requirement for this approach is that the following expected values *exist*

$$\hat{u}(\xi, t) = h(\xi) E_{\xi_\theta = \xi} \chi(\theta, t). \tag{61}$$

The simplest proof of this uses the strong Markov property to check that if the indicated expected values are finite then (FNS) is satisfied. Existence of such expected values has been obtained by Bhattacharya et al. [6] through constructions of a particular class of the Fourier multipliers, referred to as *majorizing kernels,* and defined below, from which one also obtains *uniqueness* in a naturally associated function space. Uniqueness is proven by arguments which exploit discrete parameter martingale structure in successive truncations by generation of the branching random walk product. The point is that uniqueness proofs are largely insensitive to the holding time distribution. In the original development of the theory the holding times between branchings were exponentially distributed with holding time parameters determined from the principal part (Laplacian) of the equation as above. However, recent extensions by Bhattacharya et al. [7], Orum [45], and Ossiander [47] have been obtained based on *semi-Markov cascades,* which in turn have lead to new representations of both local and global multiplicative cascade solutions to (1).

A first order approach to obtain finite expected values of the branching random walk cascade will be seen to result from the observation that the product $\otimes_\xi$ satisfies $|w \otimes_\xi z| \leq |w||z|, w, z \in C^n$, and the coefficients $m(\xi)$ may be controlled by selecting Fourier multipliers such that $m(\xi) \leq 1$ or, equivalently,

$$h * h(\xi) \leq B|\xi| h(\xi), \quad \xi \in W_h, \tag{62}$$

for suitable $B > 0$. Such $h$ is referred to as a *majorizing kernel* with constant $B$. The choice of $B$ is linked to the size of $|\hat{u}_0|$ and may be standardized to $B = 1$ in this context by consideration of suitable initial data. In any case, $h$ is a standardized majorizing kernel if and only if for some, and hence any, $B > 0$ $Bh$ is a majorizing kernel with constant $B$.

The following slightly more general definition is suitable for extensions to generalized Navier-Stokes equations with fractional Laplacian, semi-linear equations and for considerations of local solutions.

**Definition.** A positive locally integrable function $h$ on $W_h \subseteq \mathbf{R}^n \backslash \{0\}$ whose closure $\overline{W}_h$ is a semi-group and such that

  i $h$ is continuous on $W_h$

 ii $h * h > 0$ a.e. on $W_h$

iii $h * h(\xi) \leq B|\xi|^\theta h(\xi)$, for $\xi \in W_h$ and some real exponent $\theta$ and some $B > 0$,

will be referred to as an *FNS-admissible majorizing kernel with constant $B$ and exponent $\theta$*. We define $h = 0$ on $W_h^c$ and refer to $W_h$ as the *support* of $h$. Those majorizing kernels $h(\xi)$ which are defined and positive for all $\xi \neq 0$ are said to be *fully supported*.

Formulated in these terms, the results of LeJan and Sznitman [37] may be interpreted in terms of two exponent one, standardized fully supported majorizing kernels, $\pi^3/|\xi|^2$ and $\alpha e^{-\alpha|\xi|}/2\pi|\xi|$. These kernels are respectively non-integrable and integrable, with equality in (iii) of the above Definition. The former is the Fourier transform (in the sense of distribution) of the Riesz potential (convolution kernel) defining the integral operator $(-\Delta)^{-1}$, while the latter is the Fourier transform of the Bessel potential defining the integral operator $(I - \Delta)^{-1}$; c.f. Folland [19, pp. 149-154].

Apart from their role in existence, uniqueness and expected value representations valid in particular function spaces, the majorizing also play a role in tracking such structure as regularity, analyticity, support size, complexification, etc. of initial data in the evolution of solutions through subspaces corresponding to such more detailed specifications as can be reflected in their Fourier transforms. More specifically, let us define a Banach space $\mathcal{F}_{h,\gamma,T}$ with a norm that depends on a Fourier multiplier $1/h$ as the completion of the set

$$\{v \in \mathcal{S}' : \hat{v}(\xi, t) = 0, \xi \in W_h^c, |v|_{\mathcal{F}_{h,\gamma,T}} = \sup_{\substack{\xi \in W_h \\ 0 \leq t < T}} \frac{|\hat{v}(\xi, t)|}{e^{-\gamma\sqrt{t}|\xi|} h(\xi)} < \infty\} \qquad (63)$$

under the indicated norm, where $\gamma \in \{0, 1\}$ serves to conveniently index two different norms we wish to consider. Here $\mathcal{S}'$ is the space of tempered distributions on $\mathbf{R}^n$. The definition of majorizing kernel and of the function spaces $\mathcal{F}_{h,\gamma,T}$ imply that the product of distributions in $\mathcal{F}_{h,\gamma,T}$ is itself a distribution. Also, implicit to the definition of the Banach space $\mathcal{F}_{h,\gamma,T}$ is the requirement that tempered distributions belonging to this space have Fourier transforms which are functions.

In passing we wish to point out that the notion of limit space referred to in Section 1 has a natural generalization in this context. In particular if $h$ is a majorizing kernel of exponent $\theta = 1$ then $h_\lambda(\xi) = \lambda^{-2} h(\xi/\lambda)$ is also a majorizing kernel of the same exponent. Moreover if $\mathbf{u} \in \mathcal{F}_{h,\gamma,T}$ then $\mathbf{u}_\lambda \in \mathcal{F}_{h_\lambda,\gamma,T}$ and $|\mathbf{u}|_{\mathcal{F}_{h,\gamma,T}} = |\mathbf{u}_\lambda|_{\mathcal{F}_{h_\lambda,\gamma,T}}$. Thus, while an exponent one majorizing kernel $h$ such that $h = h_\lambda$ defines a limit space $\mathcal{F}_{h,\gamma,T}$ in the usual sense, the norm relation

$$|\mathbf{u}|_{\mathcal{F}_{h,\gamma,T}} = |\mathbf{u}_\lambda|_{\mathcal{F}_{h_\lambda,\gamma,T}} \qquad (64)$$

is more general.

We will conclude this section with statements of general existence/uniqueness theorems possible within majorizing spaces. The main results use majorizing kernels of different exponents to establish existence, uniqueness and regularity properties of the solutions of the (FNS). Moreover these solutions have the expected value representation in terms of the multiplicative stochastic functional $\chi(\theta, t)$ of a multitype branching random walk in Fourier wavenumber space. This will be followed with some basic structure theorems which are useful for the construction of majorizing kernels.

The first result illustrates the use of a majorizing kernel $h$ of exponent 1, to obtain existence and uniqueness of the solution given by the stochastic cascade representation for small enough initial data and forcing on a time interval that is solely constrained by the length of time for which the forcing remains small. Specifically one has the following theorem.

In the statements of these results, $(-\Delta)^{-\beta}, 0 < \beta < 3/2$ denotes the negative power of the Laplacian defined as the integral operator with symbol $|\xi|^{-2\beta}$; i.e. the convolution with a Riesz potential.

**Theorem 3.1.** *Let $h(\xi)$ be a standard majorizing kernel with exponent $\theta = 1$. Fix $0 < T \le +\infty$. Suppose that $|u_0|_{\mathcal{F}_{h,0,T}} \le (\sqrt{2\pi})^3 \nu/2$ and $|(-\Delta)^{-1}g|_{\mathcal{F}_{h,0,T}} \le (\sqrt{2\pi})^3\nu^2/4$. Then there is a unique solution $u$ in the ball $\mathcal{B}_0(0,R)$ centered at 0 of radius $R = (\sqrt{2\pi})^3\nu/2$ in the space $\mathcal{F}_{h,0,T}$. Moreover the Fourier transform of the solution is given by $\hat{u}(\xi, t) = h(\xi)E_{\xi_\theta = \xi}\chi(\theta, t), \xi \in W_h$.*

As already noted the existence question is made trivial by the stochastic iteration so long as the indicated expected values are defined. A proof of uniqueness can be based on a discrete parameter martingale construction by truncated genealogy adapted from LeJan and Sznitman [37] which is used to obtain a bound on the difference between two solutions in this space by $2P(\mathcal{T} > n)$ for all $n = 0, 1, 2, \ldots$, where $\mathcal{T}$ denotes the time to extinction of the underlying discrete parameter critical binary branching process.

It should be remarked that some regularity properties of the solutions can be inferred from the particular majorizing kernel being used via Payley-Wiener theory, Sobolev embedding, etc. For example, note that the majorizing kernel $h_0(\xi) = \pi^3/|\xi|^2$ gives existence and uniqueness, but no control over tracking regularity properties of the solutions from those of initial data. However, solutions obtained using the majorizing kernels $h_\beta^{(\alpha)} = |\xi|^{\beta-2}e^{-\alpha|\xi|^\beta}, 0 < \beta \le 1, \alpha > 0$, maintain the $C^\infty-$ regularity of the initial data. as can be seen from the bound on the Fourier transform of the solution Moreover one may construct, for example, smooth compactly supported initial data whose evolution is more efficiently dominated by a kernel $h_\beta^{(\alpha)}$ with $\beta < 1$, whose decay at infinity is consistent with the compactness of support property; see Bhattacharya and Rao [5, pp. 83-89].

On the other hand, as will be elaborated below, working in the function spaces $\mathcal{F}_{h,1,T}$ it is possible to use majorizing kernels to obtain spatial analyticity of the solution. However, it should be noted that the size constraints imposed on the initial data and forcing to obtain analyticity are substantially more restrictive than those required in Theorem 3.1; i.e. existence and uniqueness of such solutions is within a smaller ball of the function space. Specifically one has

**Theorem 3.2.** *Let $h(\xi)$ be a standard majorizing kernel with exponent $\theta = 1$. Fix $0 < T \leq +\infty$. Assume $|e^{\nu t \Delta} u_0(x)|_{\mathcal{F}_{h,1,T}} \leq \frac{(\sqrt{2\pi})^3}{2} \rho \nu e^{-1/2\nu}$ and that $|(-\Delta)^{-1} g(x,t)|_{\mathcal{F}_{h,1,T}} \leq \frac{(\sqrt{2\pi})^3}{4} \rho \nu^2 e^{-1/2\nu}$ for some $0 \leq \rho < 1$. Then there is a unique solution $u$ in the ball $\mathcal{B}_1(0,R)$ centered at $0$ of radius $R = (\rho/2)(\sqrt{2\pi})^3 \nu \times e^{-\frac{1}{2\nu}}$ in the space $\mathcal{F}_{h,1,T}$.*

Under the conditions of Theorem 3.2 the asserted solution satisfies the following decay condition

$$\sup_{0 \leq t < T} \sup_{\xi \in \mathbf{R}^3} \frac{e^{\sqrt{t}|\xi|}|\hat{u}(\xi,t)|}{h(\xi)} < \infty \tag{65}$$

Thus Theorem 3.2 provides another approach generalizing that of Lemarié-Rieusset [35], Foias and Temam [20] to obtain conditions for regularity in the stronger form of spatial analyticity. In particular, by means of the following key inequality along the lines of Foias and Temam [20] and Lemarié-Rieusset [35], for $0 \leq s \leq t, \xi, \eta \in \mathbf{R}^3$,

$$\exp\{-\nu s|\xi|^2 - \sqrt{t-s}|\xi - \eta| - \sqrt{t-s}|\eta|\} \leq e^{\frac{1}{2\nu}} e^{-\sqrt{t}|\xi|} e^{-\nu s|\xi|^2/2} \tag{66}$$

the majorization of $\hat{u}(\xi,t)$ can be shifted to $e^{\sqrt{t}|\xi|}\hat{\mathbf{u}}(\xi,t)$; see Bhattacharya et al. [6, Proposition 4.1, p. 5031]. As a consequence, for example, if $e^{-\delta|\xi|}h(\xi) \in L^1$ for some $\delta$ then, for times $\sqrt{t} > \delta$, $\mathbf{u}(\mathbf{x}+i\mathbf{y},t)$ is complex analytic in the region $\mathbf{x} + i\mathbf{y} \in \mathbf{C}^3$ such that $|\mathbf{y}| \leq \sqrt{t} - \delta$. Specifically,

$$
\begin{aligned}
\mathbf{u}(\mathbf{x}+i\mathbf{y},t) &= \int_{\mathbf{R}^3} e^{i(\mathbf{x}+i\mathbf{y})\cdot\xi} \hat{\mathbf{u}}(\xi,t) d\xi \\
&= \int_{\mathbf{R}^3} e^{i\mathbf{x}\cdot\xi} e^{-(\mathbf{y}\cdot\frac{\xi}{|\xi|}+\sqrt{t})|\xi|} e^{\sqrt{t}|\xi|} \hat{\mathbf{u}}(\xi,t) d\xi
\end{aligned} \tag{67}
$$

is a complex analytic extension since

$$
\begin{aligned}
|e^{-(\mathbf{y}\cdot\frac{\xi}{|\xi|}+\sqrt{t})|\xi|} e^{\sqrt{t}|\xi|} \hat{\mathbf{u}}(\xi,t)| &\leq e^{-(\mathbf{y}\cdot\frac{\xi}{|\xi|}+\sqrt{t}-\delta)|\xi|} e^{-\delta|\xi|} h(\xi) \\
&\leq e^{-\delta|\xi|} h(\xi)
\end{aligned} \tag{68}
$$

for $|\mathbf{y}| \leq \sqrt{t} - \delta$. In particular, in the case $e^{-\delta|\xi|} h(\xi) \in L^1$ for all $\delta > 0$, analyticity holds for all times $t > 0$. This is the case for $h(\xi) = 1/|\xi|^2$ treated by Lemarié-Rieusset [35], as well as for such majorizing kernels larger at infinity as indicated in examples below; c.f. Example (viii).

Two additional corollaries to the theory for global solutions follow immediately from the stochastic cascade representation: (a) time-asymptotic solutions, c.f. Bhattacharya et al. [7], Cannone [9] and references therein, and (b) self-similar solutions, c.f. LeJan and Sznitman [37], Cannone [9] and references therein. Specifically, for the case of (a), under the majorizing conditions of Theorem 3.1 with $T = \infty$, suppose further that $\lim_{t\to\infty} \hat{g}(\xi,t) = \hat{g}_\infty(\xi)$ exists for each $\xi \neq 0$. Then, since the underlying discrete parameter binary branching is critical, $\lim_{t\to\infty} \curlyvee(\theta,t)$ exists a.s. as a finite random product. Moreover, under the majorizing conditions one has for each $t \geq 0$, with probability one $|\chi(\theta, t)| \leq 1$. Thus, by Lebesgue's Dominated Convergence Theorem

$$\chi_\infty(\xi) := \lim_{t \to \infty} \chi(\xi, t)$$

exists for each $\xi$. Now again apply Dominated Convergence to $(FNS)_h$ to obtain

$$\chi_\infty(\xi) = \int_0^\infty \nu|\xi|^2 e^{-\nu|\xi|^2 s} \left\{ \frac{1}{2} m(\xi) \int_{\mathbf{R}^3} \chi_\infty(\eta) \otimes_\xi \chi_\infty(\xi-\eta) q_\xi(d\eta) + \frac{1}{2}\varphi_\infty(\xi) \right\} ds. \tag{69}$$

Multiplication through by $h(\xi)$ establishes that

$$\hat{u}_\infty(\xi) := \lim_{t \to \infty} \hat{u}(\xi, t)$$

exists and satisfies the steady state Navier-Stokes

$$\hat{u}_\infty(\xi) = \int_0^\infty e^{-\nu|\xi|^2 s} \left\{ |\xi|(2\pi)^{-\frac{3}{2}} \int_{\mathbf{R}^3} \hat{u}_\infty(\eta) \otimes_\xi \hat{u}_\infty(\xi-\eta) d\eta + \hat{g}_\infty(\xi) \right\} ds. \quad (FNS)_\infty.$$

With regard to (b), mild self-similar solutions to Navier-Stokes of the form

$$\mathbf{u}(\mathbf{x}, t) = \frac{1}{\sqrt{t}} \mathbf{U}(\frac{\mathbf{x}}{\sqrt{t}}), \quad p(\mathbf{x}, t) = \frac{1}{t} P(\frac{\mathbf{x}}{\sqrt{t}})$$

were originally considered by Leray [36] in connection with backward-in-time evolution of Navier-Stokes as an unrealized possible approach to determine singular flow. Cannone [9] provides a brief historical overview and ultimate demise of backward self-similarity, but the existence of forward self-similarity. In this regard the simplicity of the branching random walk cascade is overwhelming. Specifically, one observes that for majorizing kernels $h$ defined by equality in (62) satisfying the limit-space scale invariance

$$h(\xi) = \lambda^{-2} h(\xi/\lambda), \quad \lambda > 0, \xi \neq 0, \tag{70}$$

it follows that (i.) $m(\xi)$ is constant and therefore trivially invariant under scale change $\xi \to \xi/\lambda, \lambda > 0$; (ii.)$\mathbf{a} \otimes_\xi \mathbf{b}$ depends on $\xi$ only through $e_\xi \equiv \frac{\xi}{|\xi|}$; under the scale change $\xi_v \to \xi_v/\lambda, \lambda > 0$, the distribution of the holding time $S_v$ changes to that of $\lambda^2 S_v$; and, finally, the spatial transition probability $\frac{h(\eta)h(\xi-\eta)}{h*h(\xi)} d\eta$ is invariant under re-scaling of Fourier frequencies. Thus taking initial data $\hat{u}(\xi)/h(\xi)$ and forcing $\hat{g}(\xi)/h(\xi)$ (constant in time) to both be homogeneous of degree 0 one obtains that

$$\hat{u}(\xi, t) = \lambda^{-2} \hat{u}(\lambda^{-1}\xi, \lambda^2 t). \tag{71}$$

Cannone [9] notes the heuristic connection between time-asymptotics and self-similar solutions obtained by observing that since for a global solution $\mathbf{u}, p$, the space-time rescaling $\mathbf{u}^{(\lambda)}(\mathbf{x}, t) := \lambda \mathbf{u}(\lambda \mathbf{x}, \lambda^2 t)$, $p^{(\lambda)}|(\mathbf{x}) = \lambda^2 p(\lambda \mathbf{x})$, is also a solution, it follows that $\mathbf{u}^{(\infty)}(\mathbf{x}, t) := \lim_{\lambda \to \infty} \mathbf{u}^{(\lambda)}(\mathbf{x}, t)$ is self-similar whenever the indicated limit exists. In particular $\lim_{t \to \infty} t^{\frac{1}{2}} \mathbf{u}(t^{\frac{1}{2}}\mathbf{x}, t) = \mathbf{u}^{(\infty)}(\mathbf{x}, 1)$.

Finally let us record *local* existence and uniqueness from more relaxed conditions on the majorizing kernels as illustrated by the following result; c.f. Bhattacharya et al. [6].

**Theorem 3.3.** *Let $h(\xi)$ be a standard majorizing kernel with exponent $\theta < 1$. Fix $0 < T \leq +\infty, \gamma \in \{0, 1\}$. Assume $e^{\nu t \Delta} u_0(x) \in \mathcal{F}_{h,\gamma,T}$ and for some $1 \leq \beta \leq 2$, $(-\Delta)^{-\frac{\beta}{2}} g(x,t) \in \mathcal{F}_{h,\gamma,T}$ Then there is a $0 < T_* \leq T$ for which one has a unique solution $u \in \mathcal{F}_{h,\gamma,T_*}$.*

Another variation on the general approach presented here involving *semi-Markov branching*, in the sense that the holding time distributions are not required to be exponentially distributed, leads to conditions for a local existence and uniqueness theory in all dimensions.; see Bhattacharya et al. [7], Orum [45], Ossiander [47].

Some sense of the class of majorizing kernels may be derived by noting from Hölder's inequality that the set of fully supported majorizing kernels with a given exponent is a log-convex set. Also if $h(\xi)$ is a majorizing kernel then so is $ce^{a \cdot \xi} h(\xi)$ for arbitrary fixed vector $a$ and positive scalar $c$. One may also show that the only fully supported homogeneous majorizing kernels in $n \geq 3$ dimensions are those of degree $n - 1$. Along these lines the Theorem 3.4 below may be viewed as something of a "toolkit" for constructing examples.

The family of standard majorizing kernels of exponent $\theta$ on $\mathbf{R}^n$ will be denoted by

$$\mathcal{H}_{n,\theta} = \{h : W_h \to (0,\infty) : h * h(\xi) \leq |\xi|^\theta h(\xi) \quad \text{for all} \quad \xi \in W_h \subseteq \mathbf{R}^n\}.$$

**Theorem 3.4.** *[Majorizing Kernel Toolkit]*

1. *Suppose that $\{q_j : 1 \leq j \leq m\}$ satisfies $q_j > 0$, and $\sum_1^m q_j = 1$. For $h_j \in \mathcal{H}_{n,\theta_j}, j = 1, ..., m$,*

$$h(\xi) = \prod_{j=1}^m h_j^{q_j}(\xi) \in \mathcal{H}_{n, \sum_{j=1}^m q_j \theta_j}$$

   *with support $W_h = \cap_{j=1}^m W_{h_j}$.*

2. *For $h_j \in \mathcal{H}_{n,\theta}, j = 1, ..., m$,*

$$h(\xi) = \prod_{j=1}^m h_j^{q_j}(\xi) \in \mathcal{H}_{n,\theta}$$

   *with support $W_h = \cap_{j=1}^m W_{h_j}$.*

3. *Fix $n \geq 1$. Suppose that $k_1, \ldots k_m$ is a partition of $n$ and for each $j = 1, \ldots, m$, $h_j$ is in $\mathcal{H}_{k_j, \theta_j}$. Then*

$$h(\xi) = \prod_{j=1}^m h_j(\xi_j), \ \xi = (\xi_1, \ldots \xi_m) \quad \text{for } \xi_j \in \mathbf{R}^{k_j}$$

   *is in $\mathcal{H}_{n, \sum_{j=1}^m \theta_j}$ with $W_h = W_{h_1} \times ... \times W_{h_m}$.*

4. *If $A$ is a $n \times n$ invertible matrix and $h \in \mathcal{H}_{n,\theta}$, then defining $\|A\| = \sup\{|A\mathbf{x}| : |\mathbf{x}| = 1\}$,*

$$h_A(\xi) := |det A| \cdot \|A\|^{-\theta} h(A\xi) \in \mathcal{H}_{n,\theta}$$

   *with support $W_{h_A} = \{A^{-1}\xi : \xi \in W_h\}$.*

5. If $h \in \mathcal{H}_{n,\theta}$ and $\psi : \mathbf{R}^n \to [0, \infty)$ satisfies $\psi(\xi) \leq \psi(\eta) + \psi(\xi - \eta)$ for all $\eta, \xi \in W_h$, then
$$h_\psi(\xi) = e^{-\psi(\xi)} h(\xi) \in \mathcal{H}_{n,\theta}.$$

6. If $h \in \mathcal{H}_{n,\theta}$, then $e^{a \cdot \xi} h(\xi) \in \mathcal{H}_{n,\theta}$ for any fixed $a \in \mathbf{R}^n$,

7. If $h \in \mathcal{H}_{n,\theta}$ then for any pseudo-metric $\rho$ on a subset of $\mathbf{R}^3$ containing $W_h$, $e^{-a\rho(\xi_0,\xi)} h(\xi) \in \mathcal{H}_{n,\theta}$ for any $a > 0$ and $\xi_0$ fixed. In particular $e^{-a|\xi|^\beta} h(\xi) \in \mathcal{H}_{n,\theta}$ for any $a > 0$ and $0 < \beta \leq 1$

In general, existence and classification of majorizing kernels is non-trivial. For example, it can be shown that any piecewise continuous $h \in \mathcal{H}_{1,1}$ must have $W_h = (0, \infty)$ or $W_h = (-\infty, 0)$. This illustrates the tradeoff between $n$ and $\theta$; if exponent $\theta > 0$, the existence of majorizing kernels with support $\mathbf{R}^n \backslash \{0\}$ is problematic for smaller values of $n$. There are however fully supported majorizing kernels of exponent $\theta = 0$ for all $n \geq 1$. Using the toolkit Theorem 3.4 one may identify a number of majorizing kernels; see Bhattacharya et al. [6].

**Examples** [Majorizing Kernels].

   i $h_1(\xi) = \frac{1}{2\pi(1+\xi^2)} \in \mathcal{H}_{1,0}$ with $W_{h_1} = \mathbf{R}$.

   ii For $n > 1$, $h_n(\xi) = (2\pi)^{-n} \prod_{j=1}^n (1+\xi_j^2)^{-1} \in \mathcal{H}_{n,0}$ with $W_{h_n} = \mathbf{R}^n$.

   iii $h_n(\xi) = \frac{\Gamma(\frac{n+1}{2})}{2\pi^{\frac{n+1}{2}} (1+|\xi|^2)^{\frac{n+1}{2}}} \in \mathcal{H}_{n,0}$

   iv For $\xi \neq 0$, $0 \leq \beta \leq 1$, $\alpha > 0$, $h_\beta^{(\alpha)}(\xi) = |\xi|^{\beta-2} e^{-\alpha|\xi|^\beta} \in \mathcal{H}_{3,1}$ with $W = \mathbf{R}^3 \backslash \{0\}$.

   v For each $\theta \in (0,1)$, $0 \leq \beta \leq 1$, and $\alpha > 0, \xi \neq 0$, one has that $h_{\theta,\beta}^{(\alpha)}(\xi) = \frac{|\xi|^{\theta(\beta-2)} e^{-\alpha\theta|\xi|^\beta}}{(2\pi)^{3(1-\theta)} \prod_{j=1}^3 (1+\xi_j^2)^{(1-\theta)}} \in \mathcal{H}_{3,\theta}$ with $W = \mathbf{R}^3 \backslash \{0\}$.

   vi For each $\theta \in (0,1)$, $0 \leq \beta \leq 1$, and $\alpha > 0, \xi \neq 0$, one has that $h_{\theta,\beta}^{(\alpha)}(\xi) = \frac{|\xi|^{\theta(\beta-2)} e^{-\alpha\theta|\xi|^\beta}}{(1+|\xi|^2)^{2(1-\theta)}} \in \mathcal{H}_{3,\theta}$ with $W = \mathbf{R}^3 \backslash \{0\}$.

   vii For $n \geq 3$ and $(\beta, \gamma)$ with $0 \leq \beta \leq 1$ and $1 \leq \gamma \leq 1 + \beta$, one has that $h_{n,\beta,\gamma}(\xi) = \int_{t>0} t^{\frac{\gamma-n}{2}-1} e^{-t^\beta - |\xi|^2/t} dt \in \mathcal{H}_{n,1}$ with $W_h = \mathbf{R}^n \backslash \{0\}$

   viii Let $G$ be defined a.e. on the unit sphere with $G(\theta) \to \infty$ as $\theta$ approaches the points $(0,0,\pm 1), (0,\pm 1,0), (\pm 1,0,0)$, respectively. Then one has that $h_{3,1}(\xi) = \frac{1}{|\xi|^2} G(\frac{\xi}{|\xi|}) \in \mathcal{H}_{3,1}$.

The families of non-fully supported non-radial kernels illustrated in Example (viii) are *larger* kernels that permit broader existence and uniqueness results for given initial data $\mathbf{u}_0$ of (FNS). In particular the growth of $h_{3,1}$ along particular directions is much larger than $h_0(\xi) = 1/|\xi|^2$. Transforming $h_{3,1}$ via a rotation permits such growth in any direction.

The kernels in Example (vii) are closely related to the Bessel kernels of Aronszajn and Smith [2]. They can also be combined with other kernels of the Example (vii) along the lines of the previous proposition to construct kernels in $\mathcal{H}_{n,\theta}$ for $0 < \theta < 1$. One may apply the Laplace method for estimating integrals to show that the Bessel type kernels $h = h_{n,\beta,\gamma}$ are also regularizing kernels in the sense that the distributions in the corresponding function space $\mathcal{F}_{h,0,T}$ are $C^\infty-$ functions.

In view of the role of majorizing kernels in providing bounds on the Fourier transformed forcing and/or initial data, the theory contains a dual problem which is to identify classes of divergence free vector fields in physical space which are so dominated. This section is concluded with some examples in this category.

The first example is a class of solenoidal vector fields on $\mathbf{R}^3$ whose Fourier transforms are dominated by $h_{3,\beta,\gamma}(\xi)$.

**Example.** Fix $0 \le \beta \le 1$ and $1 \le \gamma \le 1+\beta$. For $1 \le j \le 3$ let $m_j(t)$ be measurable functions on $[0,\infty)$ such that $|m_j(t)| \le t^{\frac{\gamma}{2}-1} e^{-t^\beta}$ and $\int_{t>0} t^{-3/2}|m_j(t)|dt < \infty$. Let $v(x)$ be the vector field whose components $v_j(x)$ are the Laplace transforms of $m_j(t)$ evaluated at $|x|^2/4$; that is

$$v_j(x) = \int_0^\infty e^{-t|x|^2/4} m_j(t) dt$$

Let $u(x)$ be the divergence free projection of $v(x)$. Then

$$|\hat{u}(\xi)| \le c h_{3,\beta,\gamma}(\xi).$$

For the next example we consider majorization by the kernels $h_\beta^{(\alpha)}$.

**Example.** Let $\mathcal{M}$ denote the space of finite signed measures on $\mathbf{R}^3$ with total variation norm $\| \ \|$. Let $0 < \beta \le 1$ and denote the "Fourier transformed Bessel kernel" of order $\beta$ by $G_\beta(x) = (1+|x|^2)^{-\frac{1+\beta}{2}}$. Then for each $g = G_\beta * \mu, \mu \in \mathcal{M}$, one has for $\beta = 1$, $\alpha \in (0,1)$ and for $\beta \in (0,1)$, $\alpha > 0$,

$$|\hat{g}(\xi)| \le C_\beta^{(\alpha)} h_\beta^{(\alpha)}(\xi) \|\mu\|, \quad \xi \ne 0,$$

for a constant $C_\beta^{(\alpha)} > 0$. In particular, if $v \in L^1$ is a divergence free vector field then $g = G_\beta * v$ is also a divergence free vector field whose Fourier transform is dominated by $h_\beta^{(\alpha)}$.

Finally the following example uses the $h_\beta^{(\alpha)}$ majorizing kernels to provide smooth solenoidal vector fields, including some with compact support, for which the theory applies.

**Example.** Let $m_j(t), t > 0, j = 1,2,3$, be measurable functions such that $\int_0^\infty e^{-|x|^2 t}|m_j(t)|dt < \infty, x \in \mathbf{R}^3, j = 1,2,3$. Define a vector field with components $v_j, j = 1,2,3$, by

$$v_j(x) = \int_0^\infty e^{-|x|^2 t} m_j(t) dt, x \in \mathbf{R}^3.$$

Let $u$ be the divergence free projection of $v$. Then,

(i.) If $|m_j(t)| \le ct^{-\frac{1}{2}}$ then $|\hat{u}_j(\xi)| \le c' h_0^{(\alpha)}(\xi)$ for some $c' > 0, j = 1,2,3$.

(ii.) If $|m_j(t)| \le ce^{-2\alpha^2 t}$ then $|\hat{u}_j(\xi)| \le c' h_1^{(\alpha)}(\xi)$ for some $c' > 0, j = 1,2,3$.

(iii.) For arbitrary $\epsilon > 0$ there is a smooth probability density function $k_\epsilon$ supported on $[-\epsilon, \epsilon]^3$ such that

$$|\hat{k}_\epsilon(\xi)| \le c(\beta, \epsilon) \exp\{-|\epsilon\xi|^\beta\}, \xi \in \mathbf{R}^3, c(\beta, \epsilon) > 0.$$

Let $v$ be any divergence-free integrable vector field such that $|\hat{v}(\xi)| \le c|\xi|^{-2}$, $\xi \ne 0$. Then the component-wise perturbation $u = k_\epsilon * v$ is a divergence-free infinitely differentiable vector field such that $|\hat{u}_j(\xi)| \le c' h_\beta^{(\alpha)}(\xi)$, for $\alpha = \epsilon^\beta$ and some $c' > 0, j = 1, 2, 3$.

## 4. Some Concluding Remarks

We wish to conclude with pointers to two major problem areas pertaining to incompressible Navier-Stokes equations which also fall within the subject area of probability theory.

Probabilistic methods in the study of Navier-Stokes equations are not new, howeverr new ideas have emerged recently which are based on objects of broad probabilistic interest such as multi-type branching random walks. The theory to date is in a form closely tied to the contraction map methods of Kato [30], so much so that most can be obtained by either route. Bounding the stochastic cascade product by one is a bit simpler than the contraction argument only in the case that the contraction is not strict. Similarly, as illustrated in this survey, owing to the a.s. finiteness of critical branching, uniqueness, time-asymptotic limits and certain self-similarity structure are extremely simple observations within the stochastic cascade framework.

The iterates of the Picard scheme can be probabilistically identified as expectations over truncations ( by generation) of the branching random walk. In particular, the $n$−th iterate is represented by an expected value restricted to averages over trees whose evolution by time $t$ consists of only vertices $v$ such that $|v| \le n$ and such that the n-th generation vertices $|v| = n$ survive beyond time $t$, i.e. $\sum_{j=0}^{n} S_{v|j} > t$; see Proposition 4.3 of Bhattacharya et al. [6]. This and the uniqueness proof are two instances in which minimal but significant consideration is given to this underlying branching evolution.

However suitable exploitation of the genealogy of the branching random walk in ways which take fuller advantage of incompressibility has the potential to yield significant new results. Problems of this latter type are indeed at the next frontier of this research.

The Fourier representations are naturally constrained to free space and/or periodic boundary conditions. A second major problem for probabilists is the continued quest for representations in physical space which are amenable to the imposition of boundary conditions of various types.

## References

[1] Anderson, J.D. (1997): A history of aerodynamics and its impact on flying machines, Cambridge. MR1476911

[2] Aronszajn, N., K.T. Smith (1961): Theory of Bessel potentials I, *Ann. Inst. Fourier (Grenoble)* **11** 385-475 MR143935

[3] Batchelor, G.K. (1977): An introduction to fluid dynamics. New York: Cambridge University Press. MR1744638

[4] Bhattacharya, R.N., E. Waymire (1990): Stochastic processes with applications, Wiley, NY. MR1054645

[5] Bhattacharya, R.N., R.R. Rao (1976) Normal approximation and asymptotic expansions, Wiley, NY. MR436272

[6] Bhattacharya, R., L. Chen, S. Dobson, R. B. Guenther, C. Orum, M. Ossiander, E. Thomann, E. C. Waymire (2003): Majorizing kernels and stochastic cascades with applications to incompressible Navier-Stokes equations, **355** (12), 5003-5040. MR1997593

[7] Bhattacharya, R., L. Chen, R. B. Guenther, C. Orum, M. Ossiander, E. Thomann, E. C. Waymire (2004): Semi-Markov cascade representations of local solutions to 3d-incompressible Navier-Stokes, in: IMA Volumes in Mathematics and its Applications **140**, Probability and partial differential equations in modern applied mathematics, eds. J. Duan and E. C. Waymire, Springer-Verlag, NY. (in press).

[8] Busnello, B. F. Flandoli, M. Romito (2003): A probabilistic representation for the vorticity of a 3D viscous fluid and for general systems of parabolic equations, arXiv:math.PR/0306075

[9] Cannone, M. (2004): Harmonic analysis tools for solving the incompressible Navier-Stokes equations, Handbook of Math. Fluid Dynamics, **3**, eds Friedlander, S. and D. Serre, Elsevier. MR2099035

[10] Cannone, M (1995): Ondellettes, paraproduits et Navier-Stokes equations, Diderot Editeur, Paris. MR1688096

[11] Cannone, M., F. Planchon (2000): On the regularity of the bilinear form for solutions of the incompressible Navier-Stokes equations in $\mathbf{R}^3$, Rev. Mat. Iberoamericana **1** 1-16. MR1768531

[12] Chen, L., S. Dobson, R. Guenther, C. Orum, M. Ossiander, E. Waymire (2003) On Itô's complex measure condition, IMS Lecture-Notes Monographs Series, Papers in Honor of Rabi Bhattacharya,eds. K. Athreya, M. Majumdar, M. Puri, E. Waymire, **41** 65-80. MR1999415

[13] Chorin A. J. (1973): Numerical study of slightly viscous flow J. Fluid Mechanics, **57** 785-796. MR395483

[14] Chorin, A. J. (1994): Vorticity and turbulence. Springer, New York. MR1281384

[15] Chorin, A.J., J.E. Marsden (1992): A mathematical introduction to fluid mechanics, 3rd edition, Springer-Verlag, NY.

[16] Esposito, R.Marra and H.T. Yau (1996): Navier-Stokes equations for stochastic particle systems on the lattice, Commun. Math. Phys. **182** 395–456. MR1447299

[17] Fefferman, C. (2000): Existence & smoothness of the Navier-Stokes equation, The Clay Mathematics Institute, http://www.esi2.us.es/~mbilbao/claymath.htm

[18] Finn, R. (1965): Stationary solutions of the Navier-Stokes equations, Proc.

Symp. Appld. Math. **19** Amer. Math. Soc., 121-153. MR182816

[19] Folland, G. (1995): Introduction to partial differential equations, 2nd edition, Princeton, NJ. MR1357411

[20] Foias, C., R. Temam (1989): Gevrey class regularity for the solutions of the Navier-Stokes equations, J. Funct. Anal. **87** 359-369. MR1026858

[21] Fujita, H., T. Kato (1964): On the Navier-Stokes initial value problme I, Arch. Rat. Mech. Anal. **16** 269-315. MR166499

[22] Furioli, G. P.G. Lemarié-Rieusset, E. Terraneo (1997): Su l'unicite' dans $L^3(\mathbf{R}^3)$ des solutions mild des equations de Navier-Stokes, C.R. Acad. Sci. Pris, t.325, Serie 1, 1253-1256. MR1490408

[23] Galdi, G.P. (1994): *An introduction to the mathematical theory of the Navier-Stokes Equations,* Springer-Verlag, NY.

[24] Goodman, J.(1987) The convergence of the random vortex method, Comm. Pure Appl. Math., **40** 189–220. MR872384

[25] Gundy, R. F. (1989): Some martingale inequalities with applications to harmonic analysis, J. Functional Analysis **87** 212-230. MR1025887

[26] Gundy, R.F., M. L. Silverstein (1981): On a probabilistic interpretation for the Riesz transforms, in: Functional Analysis in Markov Processes, Proceedings, Katata and Kyoto, Lecture Notes in Mathematics **923** Springer-Verlag, NY. MR661625

[27] Gundy, R., M., N.Th. Varopoulos (1979): Les transformations de Riesz et les integrales stochastiques, C.R. Acad. Sci. Paris Ser. A **289** 13-16. MR545671

[28] Heywood, J. (1970): On stationary solutions of the Navier-Stokes equations as limits of non-stationary solutions, Arch. Rational Mech. Anal. **37** 48-60. MR412639

[29] Itô, K.(1965): Generalized uniform complex measures in the Hilbertian metric space with the application to the Feynman integral, Proc. Fifth Berkeley Symp. Math. Stat. Probab. II, 145-161. MR216528

[30] Kato, T. (1984): Strong $L^p$ solutions of the Navier-Stokes equations in $\mathbf{R}^m$ with applications to weak solutions, *Math. Z.* **187** 471-480. MR760047

[31] Kolokoltsov, V.N. (2000): Semiclassical analysis for diffusions and stochastic processes, Springer Lecture Notes in Mathematics, **1724**, Springer-Verlag, NY. MR1755149

[32] Ladyzhenskaya, O.A. (1969): The mathematical theory of viscous incompressible flows, 2nd ed., Gordon and Breach, NY. MR254401

[33] Ladyzhenskaya, O.A. (2003): Sixth problem of the millennium: Navier-Stokes equations, existence and smoothness, Uspekhi Mat. Nauk. **58**:2, 45-78. MR1992564

[34] Landau, L.D. and Lifschitz, E.M. (1987): Fluid mechanics, 2nd ed. Course of Theoretical Physics Vol **6**.

[35] Lemarié-Rieusett, P.G. (2000): Une remarque sur l'analyticité des soutions milds des équations de Navier-Stokes dans $\mathbf{R}^3$, *C.R. Acad. Sci. Paris,* t.330, Série 1, 183-186. MR1748305

[36] Leray, J. (1934): Su le mouvement d'un liquide visqueux empissant l'espace, Acta. Math. **63** 193-248.

[37] LeJan, Y. and A.S. Sznitman (1997): Stochastic cascades and 3-dimensional Navier-Stokes equations, *Prob. Theory and Rel. Fields* **109** 343-366.. MR1481125

[38] Long, D.G. (1988): Convergence of the random vortex method in two dimensions J. Amer. Math. Soc, **1** 779-804. MR958446

[39] Marchioro, C., M. Pulverenti (1989): Hydrodynamics in two dimensions and vortex theory, Comm. Math. Phys. **84** 483-503. MR667756

[40] McKean, H.P. (1975): Application of Brownian motion to the equation of Kolmogorov-Petrovskii-Piskunov, Comm.Pure. Appl.Math., **28** 323-331. MR400428

[41] Meleard, S. (2000): A trajectorial proof of the vortex method for the two-dimensional Navier-Stokes equation Ann. Appl. Probab. **10**(4), 11971211.

[42] Meyer, Y. (1993): Wavelets and operators, Cambridge Studies in Adv. Math, **37** Cambridge University Press, Cambridge.

[43] Meyer, Y. (2004): Oscillating patterns in some nonlinear evolution equations, in *Mathematical foundations of turbulent viscous flows*, eds. M. Cannone, T. Miiyakwaw, Springer Lecture Notes in Mathematics, Springer-Verlag, NY. MR1852741

[44] Montgomery-Smith, S. (2001): Finite time blow up for a Navier-Stokes like equation, *Proc. A.M.S.* **129**, 3017-3023. MR1840108

[45] Orum, C. (2004): Branching processes and partial differential equations, PhD Thesis, Oregon State University.

[46] Oseen, F.K.G. (1927): Neuere methoden und ergebnisse in der hydrodynamik, Academische Velagsgesellschaft, Leipzig.

[47] Ossiander, M. (2005): A probabilistic representation of solutions of the incompressible Navier-Stokes equations in $\mathbf{R}^3$, Prob. Theory and Rel. Fields, (to appear).

[48] Peetre, J. (1976): New thoughts on Besov spaces, Duke University Mathematics Series I, Durham NC. MR461123

[49] Ramirez, J. (2004): Monte-Carlo imulation of multiplicative cascades (application to two PDEs), MS Thesis, Oregon State University.

[50] Solonnikov (1964): Estimates for solutions of a non-stationary linearized system of Navier-Stokes equations, Trudy Math. Inst. Steklov, **70** 213-317. MR171094

[51] J. Szumbarski, P. Wald (1996): The stochastic vortex simulation of an unsteady viscous flow in a multiconnected domain, ESIAM Proc., **1** 153-167. MR1443220

[52] Temam, R. (1995): *Navier Stokes equations and nonlinear functional analysis,* SIAM, Philadelphia, PA. MR1318914

[53] Thomann, E., R. Guenther (2004): The fundamental solution of the linearized Navier-Stokes equations for spinning bodies in three spatial dimensions - time dependent case, J. Math. Fluid Mechanics, (in press).

[54] Woyczynski,W., P. Biler, and T. Funaki (1998): Fractal Burgers equations, *J. Diff. Equations* **148**, 9-46. MR1637513

[55] Yau, H.T. (1999): Scaling limit of particle systems, incompressible Navier-Stokes equations and Boltzmann equation. In:Proceedings of the Inter-

national Congress of Mathematics, Berlin 1998, **3** 193-205, Birkhauser. MR1648154

[56] Waymire, E. (2002): Lectures on multiscale and multiplicative processes in fluid flows, MaPhySto Lecture Notes NO 11, Aarhus University, Aarhus, Denmark.